\documentclass[10pt,a4paper]{amsart}

\usepackage[usenames,dvipsnames,svgnames,table]{xcolor}
\usepackage[latin1]{inputenc}
\usepackage{amsmath}
\usepackage{amscd}
\usepackage{amsfonts}
\usepackage{amssymb}
\usepackage{amsthm}
\usepackage[colorlinks=true]{hyperref}
\usepackage{graphicx} 
 
\hypersetup{urlcolor=blue, citecolor=red}

\definecolor{darkred}{rgb}{0.8,0,0.2}

\newtheorem{theorem}{Theorem}[section]
\newtheorem{corollary}{Corollary}[section]

\newtheorem{lemma}{Lemma}[section]
\newtheorem{proposition}{Proposition}[section]

\theoremstyle{definition}
\newtheorem{definition}[theorem]{Definition}

\newcommand{\basinm}{\mathcal{A}_{\infty,m}}

\newcommand{\Jm}{\mathcal{J}_m}
\newcommand{\C}{\mathbb{C}}

\def \Pm {\{P_m \}_{m=1}^\infty}

\def \Pmt {\{\tilde P_m \}_{m=1}^\infty}

\newcommand{\D}{{\mathbb D}}

\def \epsilon {\varepsilon}

\def \phi {\varphi}

\def \Jm {{\mathcal J}_m}

\def \Jmn {{\mathcal J}_m^n}

\def \Km {{\mathcal K}_m}

\def \Am {{\mathcal A}_{\infty, m}}

\def \C {\mathbb C}

\def \cbar {{\overline {\mathbb C}}}

\def \R {\mathbb R}

\def \N {\mathbb N}

\def \D {\mathbb D}

\def \lcd {l^\infty (\C^{d+1})}

\def \Pm {\{P_m \}_{m=1}^\infty}

\def \Pmt {\{\tilde P_m \}_{m=1}^\infty}

\def \Pmn {\{P_m^n \}_{m=1}^\infty}

\def \phim {\{\phi_m\}_{m=0}^\infty}

\def \C {\mathbb C}
\def \D {\mathbb D}

\def \R {\mathbb R}

\def \Z {\mathbb Z}
\def \im  {{\mathbf i}}

\def \dz {{\mathrm d}z}

\title{Orbit Portraits in Non-autonomous Iteration}
\author{Mark Comerford and Todd Woodard}
\email{mcomerford@math.uri.edu} \email{tw.zero1oo@gmail.com }



\begin{document}
\maketitle

\begin{abstract}We extend the definition of an orbit portrait to the context of non-autonomous iteration, both for the combinatorial version involving collections of angles and for the dynamic version involving external rays where combinatorial portraits can be realized by the dynamics associated with sequences of polynomials with suitably uniformly bounded degrees and coefficients. We show that, in the case of  sequences of polynomials of constant degree, the portraits which arise are eventually periodic which is somewhat similar to the classical theory of polynomial iteration. However, if the degrees of the polynomials in the sequence are allowed to vary, one can obtain portraits with complementary arcs of irrational length which are fundamentally different from the classical ones.
\end{abstract}

\section{Introduction}  

\subsection{Background}

The use of angles of external rays to give a combinatorial description of parameter spaces associated with the iteration of various complex analytic functions, most particularly quadratic polynomials, is a standard tool in complex dynamics e.g. \cite{DH, Mil2}. Sester in \cite{Ses2} was the first person to consider the matter of external rays for some form of non-autonomous polynomial iteration, specifically fibered quadratic polynomials. Some of these results were extended by the authors in \cite{CW} where it was shown using holomorphic motions how the landing points of dynamic rays 
associated with the non-autonomous iteration of a sequence of polynomials with suitably bounded degrees and coefficients
move holomorphically in a neighbourhood of a hyperbolic parameter value and how the angles of these rays are also locally preserved. 

In this paper we consider the extension to non-autonomous polynomial iteration of the concept of an orbit portrait, both as a combinatorial object involving collections of angles in the unit circle and as a dynamical object involving external rays. Taken with our previous work, this represents the first steps in a program aimed at understanding the parameter spaces for non-autonomous iteration and we intend to advance this program further in subsequent papers. 

In Theorem 4.1 of \cite{CW} we stated how, for a hyperbolic sequence of polynomials with connected iterated Julia sets, the angles of rays landing at a given point in the Julia set will be preserved on an entire hyperbolic component in parameter space. This raised the possibility that the combinatorics associated with angles of external rays remains a valid tool for describing parameter spaces for non-autonomous polynomial iteration, even though these spaces are generally infinite-dimensional. 

We begin by reminding the reader of the basic definitions for non-autonomous polynomial iteration and then recall certain key concepts such as hyperbolicity  which will be of importance in this context. Next we define what an orbit portrait is in the non-autonomous context, both in the combinatorial and dynamic settings. Having done this, we then state the two main theorems of this paper mentioned in the abstract: firstly that, in the case of constant degree, all non-autonomous portraits are eventually periodic which is analagous but, as we shall see, not identical to the classical case and secondly that, in the case of varying degree, one can find portraits with irrational complementary angle arcs which is impossible in the classical case. The rest of the paper (Section 2) is then devoted to proving these results.

\subsection{Non-Autonomous Iteration}

Let $d \ge 2$, $K \ge 1$, $M \ge 0$ and let $\Pm$ be a sequence of polynomials where each
\vspace{.2cm}
\[P_m(z) = a_{d_m,m}z^{d_m} + a_{d_m-1,m}z^{d_m-1} + \cdots \cdots +
a_{1,m}z + a_{0,m}\vspace{-.2cm}\]

is a polynomial of degree $2 \le d_m \le d$ whose coefficients satisfy
\vspace{.2cm}
\[\qquad 1/K \le |a_{d_m,m}| \le K,\quad m \ge 1, \:\:\quad \quad |a_{k,m}| \le M,\quad m \ge 1,\:\: 0 \le k \le d_m -1. \]

\vspace{-.2cm}
We call such sequences \emph{bounded polynomial sequences} or simply \emph{bounded sequences} and we will refer to the numbers $d$, $K$, $M$ as the \emph{bounds} for the sequence $\Pm$. 

\vspace{0cm}
For $1 \le m$ denote by $Q_m$ the composition $P_m \circ \cdots \cdots \circ P_2 \circ P_1$ and for $1 \le m \le n$ by 
$Q_{m,n}(z)$ the composition $P_{n} \circ \cdots \cdots \circ P_{m+2} \circ P_{m+1}$, (where $Q_{m,m}$ is simply the identity). Let $D_m$ and $D_{m,n}$ denote the degrees of $Q_m$ and $Q_{m,n}$ respectively, so that $D_m = \prod_{i=1}^m d_i$ and $D_{m,n} = \prod_{i=m+1}^n d_i$. If $\Pm$ is a bounded sequence, it is easy to see that we can find $R > 0$ such that, for all $m \ge 0$, if $|z| > R$, then $|Q_{m,n}(z)| \to \infty$ as $n \to \infty$. Such a radius is the called \emph{an escape radius} for the sequence $\Pm$. Note that we can find an escape radius $R$ which depends only on the bounds $d$, $K$, $M$ for our sequence and which works for every sequence which satisfies these bounds. For each $m \ge 0$ we can then define the sets
\vspace{.2cm}
\begin{eqnarray*}
\Km & = & \{z \in \mathbb{C}: \limsup_{n \to \infty}|Q_{m,n}(z)| < \infty\}, \\
\Am & = & \{z \in \mathbb{C} : \lim_{n \to \infty} |Q_{m,n}(z)| = \infty\},\\
\Jm & = & \partial \mathcal{K}_m = \partial \basinm.
\end{eqnarray*}

\vspace{-.2cm}
Here ${\mathcal K}_m$ is called the \emph{$m^{th}$ iterated filled Julia set}, $\mathcal{A}_{\infty,m}$ is the \emph{$m^{th}$ iterated basin of attraction of infinity}, and $\mathcal{J}_m$ is  the \emph{$m^{th}$ iterated Julia set}.  We also have the Green's function 
with pole at infinity for the corresponding iterated basin of infinity $\Am$ which we shall denote by $G_m$ (see \cite{Com4} for details and basic properties).

As in the classical theory, the \emph{$m^{th}$ iterated Fatou set} $\mathcal{F}_m=\mbox{int}(\mathcal{K}_m) \cup \basinm$ is the domain of normality for the family of functions $\{Q_{m,n}\}_{n=m+1}^\infty.$  It is easy to show that these sets are forward and backward invariant in the sense that, for any $0 \le m \le n$, $Q_{m,n}({\mathcal F}_m)={\mathcal F}_n$ and $Q_{m,n}({\mathcal J}_m)={\mathcal J}_n$ and $Q_{m,n}$ maps components of ${\mathcal F}_m$ surjectively onto components of ${\mathcal F}_n$. Additonally, if all the iterated Julia sets $\Jm$ are connected, we will refer to $\Pm$ as a \emph{connected sequence of polynomials}. 

If $\Pm$ is a bounded polynomial sequence as above, a sequence of points $\{z_m\}_{m=0}^\infty$ will be called an \emph{orbit} under $\Pm$ if it is invariant in the sense that, for each $m \ge 0$, $P_{m+1}(z_m) = z_{m+1}$. 

\subsection{Hyperbolicity}

We call a bounded sequence of polynomials $\{P_m\}_{m=1}^\infty$ \emph{hyperbolic} if it is uniformly expanding on its iterated Julia sets; that is, if there exist constants $C>0, \mu>1$ such that for all $i,m \ge 0$ and $z \in \mathcal{J}_m$,

\vspace{-.2cm}
\[|Q^\prime_{m,m+i}(z)|\geq C\mu^i.\]

\vspace{0cm}
For convenience, if $\{P_m\}_{m=1}^\infty$ is bounded and hyperbolic as above, we shall refer to the numbers $C$, $\mu$ as the \emph{hyperbolicity bounds} associated with $\{P_m\}_{m=1}^\infty$.

Recall that we say that a sequence $\{\{P_m^n\}_{m=1}^\infty\}_{n=1}^\infty$ of bounded sequences converges \emph{pointwise} to another limit sequence $\{P_m\}_{m=1}^\infty$ if there exists $2 \le d$ such that every polynomial $P_m^n$ has degree $2 \le d_{m,n} \le d$, every polynomial $P_m$ has degree $2 \le d_m \le d$ and, for each $m \ge 1$ and each $0 \le i \le d$, the associated coefficients $a_{i,m}^n$ for $P_m^i$ converge to the corresponding coefficient $a_{i,m}$ of $P_m$ (some of these coefficients being $0$ in the case of polynomials of degree less than $d$).

\vspace{.1cm}
\begin{theorem}  [\cite{Com4} Corollary 3.2] Let $\{\{P_m^n\}_{m=1}^\infty\}_{n=1}^\infty$ be a sequence of bounded polynomial sequences which are bounded and hyperbolic with the same constants $d$, $K$, $M$, $C$, $\mu$ and let $\Jmn$ be the corresponding iterated Julia sets. Suppose also that $\{P_m^n\}_{m=1}^\infty$ converges pointwise to a bounded sequence $\{P_m\}_{m=1}^\infty$ with iterated Julia sets $\Jm$. Then $\Pm$ is also hyperbolic with these constants and, for each $m \ge 0$, ${\mathcal J}^n_m \to {\mathcal J}_m$ in the Hausdorff topology as $n \to \infty$.
\end{theorem}

\vspace{-.2cm}
This also follows from a result of Sumi (\cite{Sum1} page 583 Theorem 2.14) as well as a result in the paper of Sester (\cite{Ses1} page 411, Proposition 4.1), both of whom were working in the context of polynomials fibered over a compact set. If $\Pm$ is a bounded sequence, for each $0 \le m < n$, let us denote by $C_{m,n}$ the set of critical values of $Q_{m,n}$ which is a set at time $n$. We then define the \emph{postcritical distance} $PD(\Pm)$ by 

\[ PD(\Pm) = \inf_{m \ge 0, n \ge m} \mathrm{dist}(C_{m,n} , {\mathcal J}_n)\]

where $\mathrm{dist}(\cdot,\cdot)$ is the usual Euclidean distance between sets. We will need the following result whose proof can be found in \cite{Com4}. We remind the reader that one condition is said to imply another \emph{up to constants} if the constants associated with the first condition give non-trivial bounds for those associated with the second. 

\vspace{.2cm}
\begin{theorem}  [\cite{Com4} Theorem 1.3] Let $\{P_m\}_{m=1}^{\infty}$ be a bounded sequence.  Then $\{P_m\}_{m=1}^{\infty}$ is hyperbolic if and only if $PD(\Pm) \ge \delta$ for some $\delta >0$. Furthermore, this equivalence is up to constants.
\end{theorem}

This result is the non-autonomous analogue of the result from classical complex dynamics that hyperbolicity is equivalent to the closure of the postcritical set being disjoint from the Julia set and also follows from the work of Sester (\cite{Ses1} page 395 Th\'eor\`eme 1.1). For complete proofs of these statements and a more detailed treatment of the Fatou-Julia theory in the setting of non-autonomous iteration, the reader is referred to \cite{BB, Com1, Com4}. 

\subsection{Orbit Portraits}

Having thus far reviewed existing material, our first new task in this paper is to give a suitable definition of an orbit portrait in the non-autonomous context. Recall that it was shown in \cite{CW} Theorem 3.1 that, for a bounded sequence $\Pm$ as above all of whose iterated Julia sets are connected, we have a sequence of B\"ottcher maps $\phi_m$ which conjugate $\Pm$ to the sequence $\{z^{d_m}\}_{m=1}^\infty$ in the sense that, 
for each $m \ge 0$, $\phi_m(z)$ is a conformal map between $\mathcal{A}_{\infty, m}$ and $\overline {\mathbb C}\setminus \overline {\mathbb D}$ for which $\phi_{m+1}\circ P_{m+1} \circ \phi_m^{\circ -1} = z^{d_{m+1}}$ (note that the original result is actually stated for a suitable parameter neighbourhood of the sequence in an analytic family, but we can conclude what we want either by considering the sequence as being in an analytic family of sequences all of whose coefficients are constant or by noting that the proof in Theorem 3.1 of the existence of the B\"ottcher map for one sequence does not require the use of an analytic family). Note also, that if the sequence $\Pm$ is monic (as will be the case with all the examples we consider in this paper), then the B\"ottcher mappings $\phi_m$ are uniquely defined by the requirement that the derivative at infinity be $1$ and the mappings are thus tangent to the identity there. 

An \emph{external ray} $R_{\theta; m}$ of angle $\theta \in \R/\Z$ is the inverse image of a straight ray from infinity under the B\"ottcher mapping $\phi_m$ at time $m$ of the form 

\vspace{-.2cm}
\[R_{\theta; m} = \phi_m^{\circ -1}(\{R e^{2\pi \im \theta}, R > 1\}).\]

The external ray $R_{\theta; m}$ is said to \emph{land at a point $p_m \in \Jm$} if 

\vspace{-.2cm}
\[ \lim_{R \to 1_+}( \phi_m^{\circ -1}(R e^{2\pi \im \theta})) = p_m.\]

Before we can define an orbit portrait associated with an orbit for a bounded polynomial sequence, we first need the definition of a formal orbit portrait involving collections of angles.
 
\vspace{.4cm}
\begin{definition} Let $2 \le d \le \infty$, let $N \ge 1$ and for each $m \ge 1$ let $d_m$ be an integer with $2 \le d_m \le d$. A \emph{formal orbit portrait $\mathcal P$ of valence $N$ and degree bound $d$} is a sequence

\vspace{-.4cm}
\[{\mathcal P} = \{A_m\}_{m=0}^\infty\]

\vspace{0cm}
where, for each $m \ge 0$, $A_m = \{\theta_{1;m},\theta_{2;m},\dots , \theta_{N;m}\}$ is an $N$-tuple of distinct angles in $\R/\Z$ and such that

\vspace{-.4cm}
\[ \theta \longmapsto d_{m+1}\theta \pmod 1\]

maps $A_m$ bijectively to $A_{m+1}$ while preserving cyclic ordering.
\end{definition}

\vspace{-.4cm}
This is an extension of the classical definition of an orbit portrait to the case where the degree of the map is allowed to vary. It is important to note that this does not directly generalize the classical definition of an orbit portrait.  Indeed, if $\{P_m\}=\{P,P,P,...\}$ for some polynomial $P$ which has a classical orbit portrait associated to some periodic orbit of period $p > 1$, then, under our definition this would yield $p$ \emph{distinct} non-autonomous portraits for the sequence $\{P,P,P,\ldots \ldots\}$. This is not entirely unexpected, as the non-autonomous definition for a grand orbit must make the same distinction when dealing with periodic orbits, e.g. Definition 2.1 on page 47 of \cite{Com1}. 

Note that the requirement for constant valence and the preservation of cyclic ordering is motivated by the same considerations as in the classical case, namely that, away from critical points, analytic functions are locally injective and sense-preserving (in our case, we will be considering portraits which are realized by hyperbolic sequences of polynomials with connected Julia sets, so these conditions will be met - see Definition 1.5 for details). 
Although our definition allows portraits of valence $1$, we will not consider these in our paper and all portraits from now on will have valence at least $2$. Since two points on the circle do not have a single cyclic ordering, the requirement for preserving cyclic ordering is vacuous in this case. Because of this, valence $2$ portraits sometimes need to be treated separately to those of higher valence as we will see later (in the proofs of Theorem 1.6 and Theorem 1.8).

We call the integers $d_{m+1}$ above the \emph{degrees} associated with the portrait ${\mathcal P}$. Note that $\theta \longmapsto d_{m+1}\theta \pmod 1$ is injective on an (open) complementary arc if and only if the length of this arc is less than or equal to $\tfrac{1}{d_{m+1}}$ (actually, there cannot be an interval of length exactly $\tfrac{1}{d_{m+1}}$ as this would map two angles of $A_m$ to the same angle of $A_{m+1}$, which violates the requirement for injectivity in Definition 1.3 above).  In the case of an arc of length $>\tfrac{1}{d_{m+1}}$ , we can find $2 \le k \le d_{m+1}$ such that the image of this arc will cover some number of adjacent complementary arcs for $A_{m+1}$ $k$ times and the other complementary arcs $k-1$ times. It follows from the fact that cyclic order must be preserved that, for the image of a single complementary arc of length $>\tfrac{1}{d_{m+1}}$, only \emph{one} arc can be covered $k$ times (note that this is true even in the case of valence $2$). This motivates the following definition.

\vspace{-.1cm}
\begin{definition}
Let ${\mathcal P} = \{A_m\}_{m=0}^\infty$ be an orbit portrait as above. For $m \ge 0$, we say an (open) complementary arc $I_m$ for $A_m$ is a \emph{critical arc} if $I_m$ has length $>\tfrac{1}{d_{m+1}}$ and the mapping $\theta \longmapsto d_{m+1}\theta \pmod 1$ is not injective on this arc. 

If $2 \le k \le d_{m+1}$ is the maximum topological degree of the restriction to $I_m$ of $\theta \longmapsto d_{m+1}\theta \pmod 1$, then the unique complementary arc $I_{m+1}$ of $A_{m+1}$ which is covered $k$ times by the image of $I_m$ is called the \emph{critical value arc for $A_{m+1}$ associated with $I_m$.}

Finally, we say ${\mathcal P}$ is \emph{unicritical} if, for each $m \ge 0$, there is just one critical arc and all but one of the complementary arcs of $A_m$ are mapped bijectively onto a (unique) complementary arc of $A_{m+1}$. 

\end{definition}

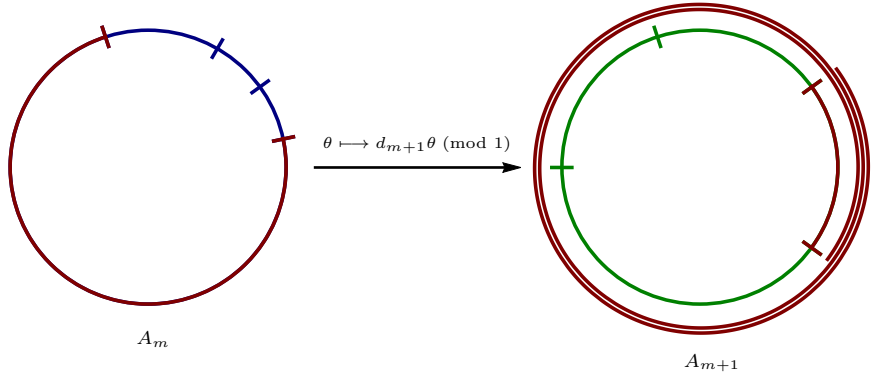
\begin{figure}[htbp]

\scalebox{0.35}{\includegraphics{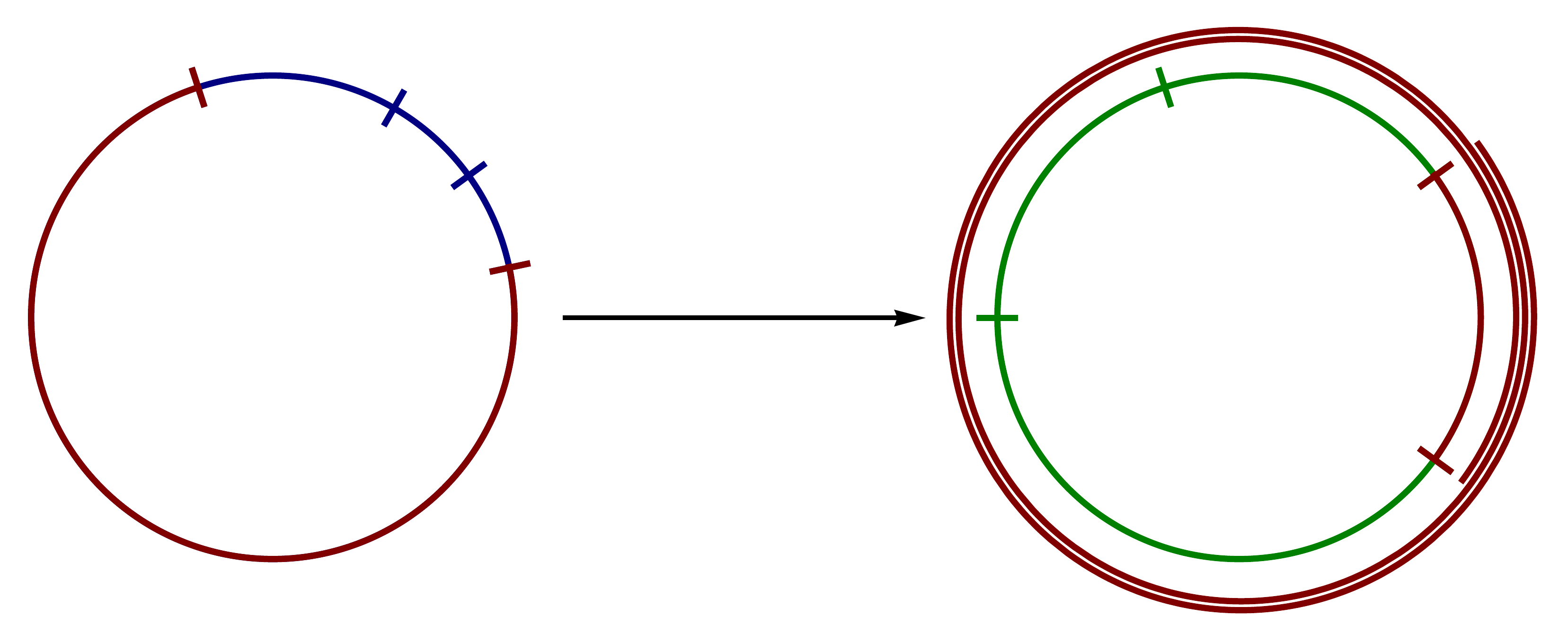}}
\unitlength1cm
\begin{picture}(0.01,0.01)
\put(-7.55,2.7){\tiny$\theta \longmapsto d_{m+1}\theta\! \pmod 1$}
\put(-10,.1){\scriptsize $A_m$}
\put(-2.8,-.2){\scriptsize $A_{m+1}$}
\end{picture}

\caption{A critical arc covering a critical value arc}
\end{figure}

Following Milnor, for a unicritical portrait ${\mathcal P}$, for a given time $m$ this allows us to speak of \emph{the} critical arc (for $m \ge 0$) and \emph{the} critical value arc (for $m \ge 1$) at a given time. 

Note that a complementary arc $I_{m+1}$ for $A_{m+1}$ as above can be associated with at most one critical arc $I_m$ as otherwise this would violate injectivity and the preservation of cyclic order for $\theta \longmapsto d_{m+1}\theta \pmod 1$ and also (in the case of valence $2$), the fact that this mapping must cover the whole of the unit circle $d_{m+1}$ times. We immediately have the following lemma. The proof is similar to that of Lemma 2.5 in \cite{Mil2}.

\vspace{.2cm}
\begin{lemma}
For any portrait ${\mathcal P} = \{A_m\}_{m=0}^\infty$ as above, the non-critical arcs of $A_m$ are mapped diffeomorphically by $\theta \mapsto d_{m+1}\theta \pmod 1$ onto distinct complementary arcs of $A_{m+1}$ and the total length of these non-critical arcs must be strictly less than $\tfrac{1}{d_{m+1}}$. On the other hand, there 
must be at least one critical arc at every time $m \ge 0$ and at least one critical value arc at every time $m \ge 1$. In the unicritical case, the critical arc has length strictly greater than $1 - \tfrac{1}{d_{m+1}}$. Its image under $\theta \mapsto d_{m+1}\theta \pmod 1$ covers one complementary arc for $A_{m+1}$ $d_{m+1}$ times and every other complementary arc for $A_{m+1}$ $d_{m+1} - 1$ times. 
\end{lemma}

{\bf Proof \hspace{.4cm}} $\theta \longmapsto d_{m+1}\theta \pmod 1$ is injective on a closed interval if and only if the length is less than $\tfrac{1}{d_{m+1}}$ and maps such an interval to an interval of $d_{m+1}$ times the length. Since cyclic order is preserved, any such complementary interval will be mapped to a single complementary interval for $A_{m+1}$ while two different such intervals will be mapped to two different intervals for $A_{m+1}$. The bound on the total length of the non-critical arcs for $A_m$ then follows immediately from this. Now this mapping wraps the unit circle $d_{m+1}$ times round itself and maps $A_m$ bijectively to $A_{m+1}$ while preserving cyclic order. Thus, by keeping track of the complementary arcs of $A_m$ and their images under $\theta \longmapsto d_{m+1}\theta \pmod 1$ as one goes around the unit circle, it is not hard to deduce that at every time $m \ge 0$ there must be at least one complementary arc on which this mapping is not injective, i.e. a critical arc. From the discussion above, this leads to the existence of an associated critical value arc at time $m+1$. $\Box$

Unicritical portraits are similar to those for classical quadratic or, more generally, unicritical polynomials in that there is always a single critical arc and a single critical value arc.
Requiring that a portrait ${\mathcal P}$ be unicritical rules out possibilities such as $\{0, \tfrac{1}{2}\}$ which is invariant under $\theta \mapsto 3\theta \,(\bmod 1)$ where we have two complementary intervals each of which covers itself just twice under iteration. The following result says that this can happen, even in the classical case. 

\vspace{.2cm}
\begin{proposition}
The cubic polynomial $P(z) = z^3 + \tfrac{3}{2}z$ has a repelling fixed point at $0$ with associated orbit portrait $\{0, \tfrac{1}{2}\}$. 
\end{proposition}

\vspace{0cm}

We remark that we cannot in general speak of a characteristic arc in the non-autonomous context as the concept of periodicity is no longer relevant. It is not even true in general that the critical value arc must be the smallest. For example, consider the quadratic portrait ${\mathcal P} =  \{A_m\}_{m=0}^\infty$ where $A_0 = \{\tfrac{1}{14}, \tfrac{1}{7}, \tfrac{2}{7}\}$ and $A_m = \{\tfrac{1}{7}, \tfrac{2}{7}, \tfrac{4}{7}\}$ for $m \ge 1$ and which is invariant under the angle doubling map. In this case the critical value arc for $A_1$ has been chosen to be $(\tfrac{4}{7}, \tfrac{1}{7})$ which is actually the largest complementary arc in $A_1$. 

As always in dynamics, the question is whether or not the combinatorial objects we introduced can be realized by a dynamical system, and so we make the following definition.

\vspace{.3cm}
\begin{definition}
A formal orbit portrait ${\mathcal P} = \{A_m\}_{m=0}^\infty$ as above is said to be \emph{realized} or a \emph{dynamic orbit portrait} if there exists a bounded polynomial sequence $\{P_m\}_{m=1}^\infty$ such that the following hold:

\begin{enumerate}
\item $P_m$ has degree $d_m$ for each $m \ge 1$,

\vspace{.2cm}
\item The iterated Julia sets $\Jm$, $m \ge 0$ are all connected, 

\vspace{.2cm}
\item There exists an orbit $\{p_m\}_{m=0}^\infty$ where, for each $m \ge 0$, $p_m \in \Jm$ and such that    
the external rays associated with the formal portrait (at time $m$), denoted

\[\{R_{\theta_{1; m}}, R_{\theta_{2; m}}, \dots \dots , R_{\theta_{N; m}}\}\]

\vspace{.3cm}
\hspace{-.4cm}(and no others) all land on $p_m$. 

\end{enumerate}

\vspace{-.2cm}
We say that an orbit portrait is a \emph{hyperbolic orbit portrait} or \emph{hyperbolically realizeable} if it is realizeable by a hyperbolic bounded connected sequence of polynomials.

\end{definition}

\vspace{-.2cm}
Note that if a portrait ${\mathcal P}$ is realized, then there can be no critical points on the orbit  $\{p_m\}_{m=0}^\infty$ above, something which is automatic in the case of a hyperbolic bounded connected sequence.

\vspace{0cm}
Since we are restricting ourselves to connected polynomial sequences, all critical points and their iterates will belong to the corresponding filled Julia sets. Hence, if a particular portrait $\mathcal P$ as above is realized by a sequence $\Pm$, for each $m \ge 0$ the critical points of $P_{m+1}$ and the critical values of $Q_m$ will avoid the rays $R_{\theta_{1; m}}, R_{\theta_{2; m}}, \dots \dots , R_{\theta_{N; m}}$. Furthermore, the fact that number of rays in the portrait at time $m$ does not depend on $m$ ensures that none of these rays can land on a critical point of $P_{m+1}$. 

For each $m \ge 0$, the rays meeting at $p_m$ divide the complex plane into $N$ simply connected regions which (following Milnor), we will refer to as the \emph{sectors} associated with the portrait at time $m$. The correpsondence between sectors and external angles is clear in view of the following lemma.

\vspace{.2cm}
\begin{lemma}
Let $\Pm$ be a bounded polynomial sequence with connected iterated Julia sets, let $m \ge 0$ and let $\theta_1$, $\theta_2$ be two angles for which the corresponding external rays $R_{\theta_1;m}$, $R_{\theta_2;m}$ land at the same point of the iterated Julia set $\Jm$ at time $m$ and divide the complex plane into two simply connected domains $D_1$, $D_2$. Then, on relabelling $D_1$, $D_2$ if needed, we may assume all points in the corresponding iterated basin of infinity $\Am$ whose external angles are in the range $(\theta_1, \theta_2)$ lie in $D_1$ while all points in $\Am$ whose external angles are in the range $(\theta_2, \theta_1)$ lie in $D_2$.
\end{lemma}

{\bf Proof \hspace{.4cm}} Since $\Jm$ is connected, for any $h >0$, the equipotential curve on which $G_m(z) = h$ is a simple closed curve. Using the B\"ottcher map $\phi_m$, we see that the external rays $R_{\theta_1;m}$, $R_{\theta_2;m}$ divide this curve into two arcs. If one of these arcs contained a point whose external angle lay in 
$(\theta_1, \theta_2)$ and another point whose external angle lay in $(\theta_2, \theta_1)$, it would follow (using the continuity of the inverse B\"ottcher map) that we could obtain a separation of this arc using these complementary intervals. The desired conclusion then follows easily. $\Box$

\vspace{.4cm}
Using the above lemma, there is a one-to-one correspondence between the sectors $S_m$ at each time $m \ge 0$ and the  
complementary intervals $I_m$ for the associated collections of angles $A_m$ for the formal portrait $\mathcal P$ at this time. For a sector $S_m$, let us call the corresponding interval $I_m$ as above the \emph{complementary interval associated with $S_m$}
 
From above, for each $m \ge 0$, both the critical point of $P_{m+1}$ and the critical values of $Q_m$ must be situated in the interiors of these sectors. Following Milnor, for $m \ge 0$, we say a sector $S_m$ at time $m$ is a \emph{critical sector} if it contains a critical point of $P_{m+1}$ while a sector $S_{m+1}$ at time ${m+1}$ is a \emph{critical value sector} if it contains a critical value of $P_{m+1}$. This terminology is justified in view of the following result which extends the first three parts of Lemma 2.11 in \cite{Mil2} (note that the fourth part of Milnor's result concerns the characteristic arc which, as we already remarked, is not relevant in the context of non-autonomous iteration).

\vspace{.2cm}
\begin{theorem}
Let $\mathcal P = \{A_m\}_{m=0}^\infty$ be a formal orbit portrait as above which is realized by a connected bounded polynomial sequence $\Pm$ with associated orbit $\{p_m\}_{m=0}^\infty$ as above. Then, for each $m \ge 0$, the following hold:

\begin{enumerate}
\item A sector $S_m$ is a critical sector if and only if the associated complementary arc $I_m$ is a critical arc. Thus, any critical sector at time $m$ has angular width strictly greater than $\tfrac{1}{d_{m+1}}$ while all non-critical sectors have combined width strictly less than $\tfrac{1}{d_{m+1}}$.

\vspace{.2cm}
\item The polynomial $P_{m+1}$ maps a small neighbourhood of $p_m$ diffeomorphically to a small neighbourhood of $p_{m+1} = P_{m+1}(p_m)$, carrying each sector based at $p_m$ locally onto a sector based at $p_{m+1}$ and preserving the cyclic order of these sectors around their base point. 

\vspace{.2cm}
\item If $S_m$ is a non-critical sector at $p_m$, then $P_{m+1}$ maps $S_m$ homeomorphically onto a sector $S_{m+1}$ at $p_{m+1}$. If $S_m$ is a critical sector, then $P_{m+1}$ maps $S_m$ onto the whole complex plane. 

\vspace{.2cm}
\item If $S_m$ is a critical sector with associated (critical) complementary arc $I_m$, then the critical value arc $I_{m+1}$ for $I_m$ is associated with a critical value sector $S_{m+1}$ which contains a critical value $v_{m+1}$ which is the image under $P_{m+1}$ of a critical point $c_m$ which lies inside $S_m$. Further, if  $2 \le k \le d_{m+1}$ is such that the image of $I_m$ under $\theta \longmapsto d_{m+1}\theta \pmod 1$ covers $I_{m+1}$ $k$ times, then the image of $S_m$ covers $S_{m+1}$ $k$ times and the rest of the complex plane $k - 1$ times.

\end{enumerate}

\end{theorem}

We draw the reader's attention to the care we took in distinguishing between a critical value sector and a critical value sector corresponding to a critical value arc for some critical arc. The issue here is that if a critical value sector $S_{m+1}$ contains a critical value $v_{m+1}$ which is the image of a critical point $c_m$ of the polynomial $P_{m+1}$, then the component of $P_{m+1}^{\circ -1}(S_{m+1})$ which contains $c_m$ may not meet the point $p_m$ on the orbit which realizes the portrait and instead meets another preimage of $p_{m+1}$. The following shows that this can indeed happen. 

\vspace{.2cm}
\begin{proposition}
There exists a non-autonomous portrait associated with a single cubic polynomial for which there is a critical value sector at time $1$ which is not associated with any critical value arc at time $1$ (or critical arc at time $0$).
\end{proposition}

\vspace{0cm}

In order to use orbit portraits to distinguish between polynomial sequences which have different behaviour, it will be convenient to equate two systems which possess orbits with similar dynamic orbit portraits, including the classical case where the portrait is associated with a periodic orbit of period two or greater.

\vspace{.2cm}
\begin{definition} 
We say that two portraits $\mathcal{P}^1 = \{\{\theta_{1;m}^1,\theta_{2;m}^1,\dots , \theta_{N;m}^1\}\}_{m=0}^\infty$ and $\mathcal{P}^2 = \{\{\theta_{1;m}^2,\theta_{2;m}^2,\dots , \theta_{N;m}^2\}\}_{m=0}^\infty$ are \emph{equivalent} if they have the same valence $N$, their respective degrees 
$d^1_{m+1}$, $d^2_{m+1}$ are equal for every $m \ge 0$ and there exist $\theta \in \R/\Z$ and non-negative integers $m_1$ and $m_2$ such that (modulo $1$) we have 

\vspace{-.4cm}
\[\{\theta_{1;m_1}^1 + \theta,\theta_{2;m_1}^1 + \theta, \dots , \theta_{N;m_1}^1 + \theta\} = \{\theta_{2;m_2}^2,\theta_{2;m_2}^2,\dots , \theta_{N;m_2}^2\}.\]

\end{definition}
In other words, $\mathcal{P}^1$ is equivalent to $\mathcal{P}^2$ if a shift of $\mathcal{P}^1$ is conjugate by rotation to a shift of $\mathcal{P}^2$.  It is straightforward to check that this equivalence does in fact define a bona fide equivalence relation on orbit portraits and, in the case of a classical system with an orbit portrait associated to a period $p$ orbit, the $p$ distinct non-autonomous portraits will all be equivalent under this relation.

It is well known in the classical (autonomous) quadratic case that every formal orbit portrait is realized by an orbit for a quadratic polynomial (and in fact a large class of quadratic polynomials - see Theorem 2.4 in \cite{Mil2}).  In the present case, if we restrict ourselves to non-autonomous bounded sequences of constant degree, then we arrive at a somewhat surprising analagous result.  

\vspace{.2cm}
\begin{theorem}  
A hyperbolically realizeable orbit portrait of constant degree $d \ge 2$ is preperiodic, i.e. it is equivalent to a periodic orbit portrait. 
\end{theorem}

\vspace{-.2cm}
Despite the above, hyperbolically realizeable portraits of constant degree are not simply the same as classical portraits. Recall that classical orbit portraits must be pairwise unlinked; that is, the collections of angles of the rays for any two different points on the orbit must lie in disjoint sub-intervals of $\R/\Z$ (see Lemma 2.3 in \cite{Mil2}). In the non-autonomous case, this need no longer be true.

\vspace{.2cm}
\begin{theorem}
There exists a non-autonomous orbit portrait which does not possess the unlinking property and which is realized by a hyperbolic connected sequence of quadratic polynomials. 
\end{theorem}

\vspace{-.3cm}
Lastly, if we relax the assumption of constant degree, then we can have portraits which look very different to the classical ones. 

\vspace{.2cm}
\begin{theorem}
There exists a bounded connected unicritical hyperbolic sequence consisting of quadratic and cubic polynomials which possesses an orbit portrait for which there is a complementary arc with irrational length. 
\end{theorem}

\section{Proofs of the Main Results}

For a simply connected domain $U$, if we let $\delta(z)$ denote the Euclidean distance to the boundary $\partial U$, recall that we have the following estimate on the hyperbolic metric $\rho(\cdot\, , \cdot)$ on $U$ in terms of $\delta(z)$ (e.g. Theorem 4.3 in \cite{CG})

\vspace{-.2cm}
\[ \frac{|\dz|}{2\delta(z)}\leq d\rho(z) \leq \frac{2|\dz|}{\delta(z)}.\]

\vspace{.3cm}

We also need the following lemma on the lengths of segments of hyperbolic geodesics, which is stated as Lemma 3.2 in \cite{CW} although the original proof goes back to the work of Hiroki Sumi in the Proof of Theorem 1.12 in \cite{Sum2} and the proof of Theorem 1.1 in \cite{Sum3}. For a bounded hyperbolic sequence $\Pm$ with connected iterated Julia sets, given $m \ge 0$
and the corresponding Green's function $G_m(z)$ with pole at $\infty$ for the iterated basin of infinity $\Am$, for a point $z \in \Am$, we denote by $\gamma_z$ the segment of the Green's line in ${\mathcal A}_{\infty, m}$ which runs from $z$ to $\partial {\mathcal A}_{\infty, m} = {\mathcal J}_m$ (and which is clearly part of an external ray). Also, for a curve $\gamma$, we denote the Euclidean arc length of $\gamma$ by $\ell(\gamma)$.  

\vspace{.2cm}
\begin{lemma} Let $\Pm$ be a bounded hyperbolic sequence all of whose iterated Julia sets are connected. Then there exist constants $C > 0$, $\alpha >0$ depending only on the degree, coefficient and hyperbolicity bounds for $\Pm$ such that, for any $m \ge0$ and $z \in {\mathcal A}_{\infty, m}$, $\ell (\gamma_z) \le C G_m(z)^\alpha$. 
\end{lemma}

\vspace{-.2cm}
As observed in \cite{CW} Corollary 3.2, it follows easily from this that for a connected hyperbolic sequence, all rays must land. One can also deduce quite easily from this estimate that all prime ends must be trivial and hence the inverse B\"ottcher maps extend continuously to $\partial (\cbar \setminus \overline \D) = {\mathrm C}(0,1)$ while the iterated basins of infinity must be locally connected (\cite{Mil1} Theorems 17.12 and 17.14). Sumi used this result in order to prove that, for a connected hyperbolic (or, in his case, more generally, semi-hyperbolic) sequence, the iterated basins of infinity were John domains (see \cite{Sum2} Theorem 1.12 and \cite{Sum3} Theorem 1.1 for details and also Theorem 4.4 in Br\"uck's paper \cite{Br} which proves this for certain sequences of quadratic polynomials). A very easy consequence of the above is the following result which we will need later in proving Theorem 1.10.

\vspace{.3cm}

\begin{corollary}
Let $\Pm$ be a connected bounded hyperbolic sequence of polynomials. Then, for every $m \ge 0$ and every $z_m \in \Jm$, at least one external ray lands at $z_m$.
\end{corollary}

\vspace{-.2cm}
{\bf Proof \hspace{.4cm}} As above, let $G_m$ be the Green's function with pole at $\infty$ for the iterated basin of infinity $\Am$.  Using the estimate on the hyperbolic metric given at the start of this section, one sees easily that $G_m$ is continuous at $z_m$ (this also follows from 
Theorem 1.4 of \cite{Com4}). Using Lemma 2.1 this then allows us to deduce the existence of a sequence of external rays whose landing points converge to $z_m$. The desired conclusion then follows from the fact that the inverse B\"ottcher maps extend continuously to the unit circle.
$\Box$

In contrast to this, we also will need the following which is another easy consequence of 
Sumi's results. This result was mentioned informally after the proof of Corollary 3.2 in \cite{CW} but, as it will be essential to us here, we now give a full statement and proof. We observe that Sumi comes quite close to stating what we need in \cite{Sum2} where he observes (in Remark 6 on page 7) that, if $V$ is a John domain, then $V$ is finitely connected at any point in $\partial V$. 

\vspace{.2cm}
\begin{lemma}
Let $\Pm$ be a connected bounded hyperbolic sequence. Then there exists $N$ which depends only on the degree, coefficient and hyperbolicity bounds for $\Pm$ such that, for any $m \ge 0$ and $z_m \in \Jm$, there are at most $N$ rays landing at $z_m$. 
\end{lemma}

\vspace{-.2cm}
{\bf Proof \hspace{.4cm}} Recall that for a sequence $\Pm$, the iterated Julia sets are all connected if and only if none of the critical points escape to infinity under iteration. Fix degree and coefficient bounds $d$, $K$, $M$ and hyperbolicity bounds $C$, $\mu$ and let $X \subset \lcd$ denote the space of all those sequences of coefficients whose corresponding polynomial sequences have these degree and coefficient bounds, are hyperbolic with bounds $C$, $\mu$ and whose iterated Julia sets are connected. It then follows easily from Theorem 1.1 that $X$ is compact with respect to the product topology inherited from $\lcd$ which is equivalent to that arising from pointwise convergence of polynomial sequences as defined earlier. 

By \cite{Sum2} Theorem 1.12 and \cite{Sum3} Theorem 1.1, we can then find $\epsilon > 0$ which depends only on the coefficient and hyperbolicity bounds for $\Pm$ such that all the iterated basins of infinity $\Am$ are $\epsilon$-John domains where we may take the centres to be $\infty$ and the John curves to be external rays. 

If $m \ge 0$ and $z_m \in \Jm$ are arbitrarily chosen, then, using B\"ottcher coordinates, one sees that the hyperbolic distance between points on any two different external rays which land at $z_m$ must tend to infinity as we approach $z_m$. Hence, using the estimate on the hyperbolic metric at the start of this section, if $r > 0$ is sufficiently small,  the intersections of these $\epsilon$ (distance) cones for these rays with the circle ${\mathrm C}(z_m, r)$ must be disjoint. Since the intersection of each such cone will take up at least angle $2\arcsin (\epsilon)$ of the angle of this circle, it follows that the number of rays which meet at $z_m$ will be bounded by $\pi/\arcsin \epsilon$ whence the result follows. $\Box$

{\bf Proof of Theorem 1.6\hspace{.4cm}} For (1), suppose a sector $S_m$ 
defined by two of the rays determined by $A_m$, the angles of the portrait ${\mathcal P}$ at time $m$, and with associated complementary arc $I_m$ contains a critical point $c_m$ of $P_{m+1}$. 
Join $c_m$ to any point $z$ in $S_m \cap \Am$ with a curve $\gamma$ in $S_m$. Clearly, since $P_{m+1}(z)$ cannot be a critical value as it lies in the iterated basin of infinity ${\mathcal A}_{\infty, m+1}$ while all the iterated Julia sets are connected, we can ensure that the only critical value of $P_{m+1}$ which meets $P_{m+1}(\gamma)$ is $P_{m+1}(c_m)$.

If $P_{m+1}(\gamma)$ meets 
the images under $P_{m+1}$ of the rays which define $\partial S_m$, then $I_m$ must be a critical arc as, in view of Lemma 1.1, the image under $\theta \longmapsto d_{m+1}\theta \pmod 1$ of a non-critical arc does not meet any of the angles in $A_{m+1}$, the angles of the portrait $\mathcal P$ at time $m+1$. Otherwise, we can assume that this does not happen and, since $P_{m+1}$ is not injective in a neighborhood of $c_m$ while the only critical value of $P_{m+1}$ which meets $P_{m+1}(\gamma)$ is $P_{m+1}(c_m)$, by analytic continuation of a suitable locally defined branch of $P_{m+1}^{\circ -1}$ from $P_{m+1}(c_m)$, we can find a different curve $\gamma'$ in $S_m$ which is a preimage of $P_{m+1}(\gamma)$ and joins $c_m$ to a point $z'$ which is a preimage of $P_{m+1}^{\circ -1}(z)$.
Since the only critical value of $P_{m+1}$ which meets $P_{m+1}(\gamma)$ is $P_{m+1}(c_m)$, it follows (e.g. from Rouch\'e's theorem) that $\gamma'$ can only meet $\gamma$ at $c_m$ and so $z'$ must be different from $z$. Here again, it follows easily from Lemma 1.1 by looking at equipotential curves that the complementary arc $I_m$ must be a critical arc.

\vspace{.2cm}
\color{black}
\begin{figure}[htbp]

\scalebox{.292}{\includegraphics{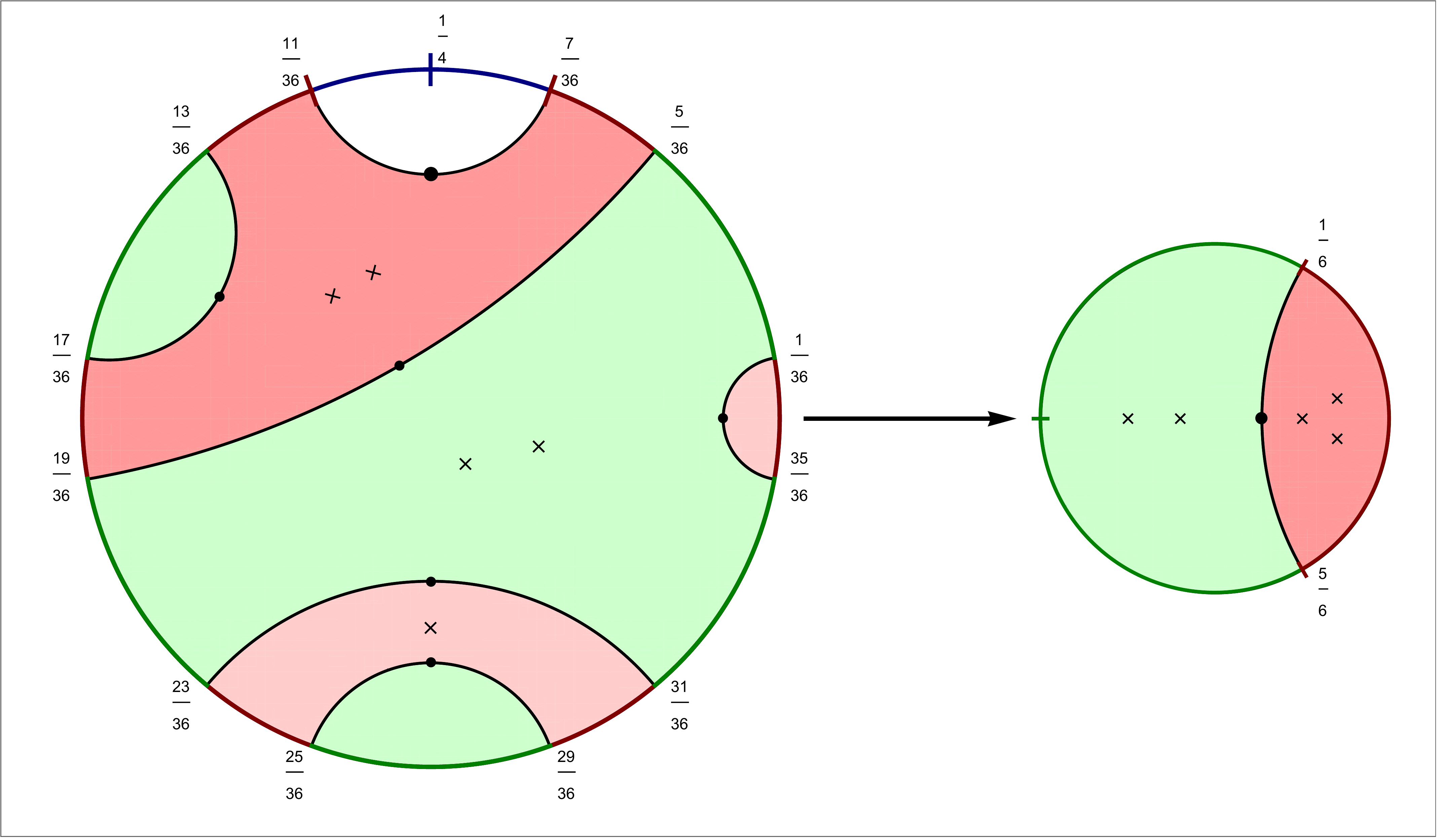}}
\unitlength1cm
\begin{picture}(0.01,0.01)
  \put(1.3,4.3){\tiny$P_{m+1}$}
  \put(-3.65,5.75){\scriptsize$U_m$}
\end{picture}

\caption{A degree $6$ example with critical arc $(\tfrac{11}{36}, \tfrac{7}{36})$ and critical value arc $(\tfrac{5}{6}, \tfrac{1}{6})$.}
\end{figure}
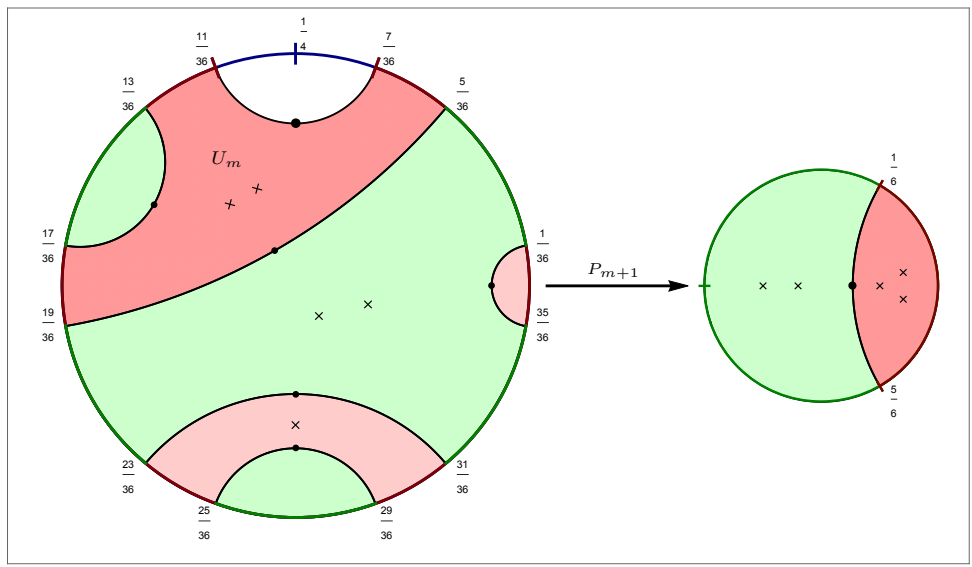

Conversely, suppose $S_m$ is a sector whose associated complementary arc $I_m$ is a critical arc (the reader may find it helpful to consult with Figure 2 below in what follows). Let $I_{m+1}$ be the (unique) critical value arc for $I_m$ and $S_{m+1}$ the corresponding sector. Since $I_m$ is a critical arc and $I_{m+1}$ is its associated critical value arc, $P_{m+1}$ maps $\partial S_m$ to $\partial S_{m+1}$, whence any component of $P_{m+1}^{\circ -1}(S_{m+1})$ must lie entirely inside or entirely outside $S_m$. Since $\theta \longmapsto d_{m+1}\theta \pmod 1$ locally preserves cyclic order, there are preimages of $I_{m+1}$ inside $I_m$ which meet each of the endpoints of $I_m$ (this can best be seen by remembering as illustrated in Figure 1 how the image of the critical arc `wraps' around the circle, covering the critical value arc one more time than the rest of the circle). Since $I_m$ is a critical arc, its image under 
$\theta \longmapsto d_{m+1}\theta \pmod 1$ covers $I_{m+1}$ at least twice and so $I_m$ contains at least two preimages of $I_{m+1}$. Since the image of the critical arc also covers the whole circle and not just the critical value arc,  it then follows that $I_m$ contains preimages of the other complementary arcs of $I_{m+1}$ and so the two preimages above which meet the endpoints of $I_m$ must be distinct. Applying the continuity of the inverse B\"ottcher mappings, we see that, for each of the rays which form part of $\partial S_m$,  
there must be at least one component of $P_{m+1}^{\circ -1}(S_{m+1})$ which is adherent to that ray. This component clearly contains points whose angles lie inside $I_m$ and so meets $S_m$. From above, we then have that it will then lie entirely inside $S_m$. On the other hand, as $p_m$ is not a critical point, the inverse image $P_{m+1}^{\circ -1}(S_{m+1})$ contains one and only one connected component which meets $p_m$.  Thus, there is one and only one preimage component which meets $p_m$, is adherent to both these rays and lies entirely inside $S_m$ (this is the darker pink component in the diagram above).

If we call this component $U_m$, then $\partial U_m$ consists of (complete) external rays and their landing points so that 
each component of the complement $\C \setminus U_m$ is clearly unbounded and thus $U_m$ must be simply connected. On the other hand, from above, the two preimages of $I_{m+1}$ which are contained in $I_m$ and meet the two endpoints of $I_m$ are distinct. 
The set of external angles of points in $U_m \cap \Am$ then includes these two preimages of $I_{m+1}$.  Thus, $P_{m+1}$ is not injective on $U_m$ and it follows from the Riemann-Hurwitz formula (e.g. \cite{Bea} page 87, Theorem 5.4.1) that $U_m$ and hence $S_m$ must contain a critical point of $P_{m+1}$. Note also that since $I_m$ contains the two distinct preimages of $I_{m+1}$ above which meet the two endpoints of $I_m$ as well as preimages of the other complementary arcs of $A_{m+1}$, it follows that $U_m$ cannot be all of $S_m$ and in particular $\partial U_m$ must contain at least one ray pair other than $\partial S_m$. 

The proof of (2) is elementary and is the same as in the classical case, depending as it does on the fact that there are no critical points on the orbit $\{p_m\}_{m=0}^\infty$.   

For (3), if we let $R_{\theta_1;m+1}$ and $R_{\theta_2;m+1}$ be the images under $P_{m + 1}$ of the two external rays which define $\partial S_m$, these rays and their common landing point $p_{m+1}$ divide the plane into two regions which we shall call $V_{m+1}$ and $W_{m+1}$. It then follows from the one to one correspondence between sectors and their associated complementary intervals and the fact that the image under $\theta \longmapsto d_{m+1}\theta \pmod 1$ of a complementary interval for $A_m$ covers all points of a given complementary interval for $A_{m+1}$ either a fixed number of times or not at all that the number of solutions in $S_m$ to the equation $P_{m+1}(z) = z_0$ will be constant as $z_0$ varies over each of $V_{m+1} \cap {\mathcal A}_{\infty, m+1}$ and $W_{m+1} \cap {\mathcal A}_{\infty, m+1}$. Since $\partial S_m$ is mapped by $P_{m+1}$ to $\partial V_{m+1} = \partial W_{m+1}$, it follows from 
Rouch\'e's theorem, local compactness and the connectedness of $V_{m+1}$ and $W_{m+1}$ that the number of solutions in $S_m$ to this equation will be constant on all of $V_{m+1}$ and $W_{m+1}$.

Since by Lemma 1.1 a non-critical arc for $A_m$ is mapped injectively to an arc for $A_{m+1}$, 
the conclusion for non-critical sectors now follows immediately from the above. 
Again by Lemma 1.1 and the above, the image of a critical sector then must include $V_{m+1} \cup W_{m+1} \cup {\mathcal A}_{\infty; m+1}$ which is all of $\C$ except for the common landing point $p_{m+1}$ of $R_{\theta_1;m+1}$ and $R_{\theta_2;m+1}$. However, if we let $U_m$ be the preimage component of $S_{m+1}$ which we used above in the proof of part (1) then, as mentioned above, $\partial U_m$ contains at least one other ray pair other than $\partial S_m$ and this ray pair must also be a preimage of $R_{\theta_1;m+1} \cup R_{\theta_2;m+1}$. The common landing point of these rays cannot also be $p_m$ as otherwise this would violate the local injectivity of $P_{m+1}$ at this point. Hence the common landing point must lie in the interior of $S_m$ and so $p_{m+1}$ is also in the range of $P_{m+1}$ when restricted to $S_m$ whence the conclusion follows. 

For the last part (4), if $S_m$ is a critical sector with associated critical arc $I_m$, if we again let $S_{m+1}$ be the sector associated with the corresponding critical value arc $I_{m+1}$ then, letting $U_m \subset S_m$ be the component of $P_{m+1}^{\circ -1}(S_m)$ in the proof of (1) above, $S_{m+1}$ contains the image of a critical point $c_m \in U_m \subset S_m$ and is thus a critical value sector with the desired properties.
When $S_m$ is a critical sector, one of the two regions $V_{m+1}$, $W_{m+1}$ in the proof of (3) above, say $V_{m+1}$, must then be the critical value sector $S_{m+1}$ corresponding to the critical value arc $I_{m+1}$. If the image of $I_m$ under $\theta \longmapsto d_{m+1}\theta \pmod 1$ covers $I_{m+1}$ $k$ times for some $2 \le k \le d_{m+1}$, then the arc $I_m$ contains $k$ preimages of $I_{m+1}$ and
$k-1$ preimages of the other complementary arcs for $A_{m+1}$. The same argument as in the proof of (3) shows that 
 the equation $P_{m+1} = z_0$ will have $k$ solutions for $z_0 \in S_{m+1}$ and $k-1$ solutions for $z_0 \in W_{m+1}$. Since the arc $I_m$ contains $k$ preimages of $I_{m+1}$, the sector 
$S_m$ must then contain $k-1$ ray pairs each of which gets mapped by $P_{m+1}$ to the boundary of the critical value sector above. Taking this together with the common landing points of these ray pairs, we see that points on the boundary of the critical value sector are covered $k-1$ times by the image of $S_m$ and the desired conclusion then follows. $\Box$

For a critical sector $S_m$ with associated critical sector $I_m$, if we let $I_{m+1}$ be the corresponding critical value arc and $S_{m+1}$ the corresponding sector as above, we just saw not only that $S_{m+1}$ is a critical value sector, but that there is one and only one component of the preimage $P_{m+1}^{\circ -1}(S_{m+1})$ (which we called $U_m$) which is adherent to all of $\partial S_m$ and must contain a critical point of $P_{m+1}$. Let us call this uniquely defined preimage component, the \emph{critical preimage component} associated with the critical sector $S_m$ (or it's corresponding critical value sector $S_{m+1}$).

Parts (3) and (4) of Theorem 1.6 imply that the image under $P_{m+1}$ of a sector $S_m$ at time $m$ is a union of sectors at time $m+1$ (or possibly the whole complex plane) and so we can say that the sectors give us a non-autonomous version of a Markov partition for our polynomial sequence. 

Furthermore, the behaviour of $P_{m+1}$ relative to the sectors is more or less completely described by diagrams as in Figure 1, one for each critical arc. These diagrams can be described much more succinctly using sets of equivalence classes, one for each critical arc, with each equivalence class consisting of a pair of angles which are preimages of the two endpoints of the critical value arc which lie in the closure of the critical arc. For example, the configuration of the preimages in the degree $6$ example in Figure 2 is specified as below: 

\[\Big \{ \big \{{\bf\{\tfrac{7}{36}, \tfrac{11}{36}\}}, \{\tfrac{13}{36}, \tfrac{17}{36}\},\{\tfrac{19}{36}, \tfrac{5}{36}\} \big \} ,\:
\big  \{ \{\tfrac{23}{36}, \tfrac{31}{36}\}, \{\tfrac{25}{36}, \tfrac{29}{36}\} \big \} ,\: \big \{\{\tfrac{35}{36}, \tfrac{1}{36}\} \big \} \Big  \}.\]

Here, the first set of pairings corresponds to the critical preimage component as defined above and the pairing which corresponds to the boundary of the critical value arc is indicated in bold. Note how these pairings must be unlinked, both within each set of such pairings and between two different sets, the reason being the same as in the classical case in \cite{Mil2}, namely that external rays are not allowed to cross. 
It is also immediate that each angle in the closure of the critical arc must belong to at most one pairing. 

The cardinality of each set of pairings is the same as the topological degree of $P_{m+1}$ on the corresponding preimage component of the critical value sector and of course subtracting one from this gives the number of critical points (counted by multiplicity). For example, in Figure 2 we see that the critical preimage component corresponds to a set with three pairings and the degree of $P_{m+1}$ on this preimage component is $3$. There are also two other sets of pairings of cardinality two and one which correspond to preimage components on which the degree of $P_{m+1}$ is $2$ and $1$ respectively. In fact, the pairings above not only determine the degrees of $P_{m+1}$ on the components of the preimages of the critical value sector, but also on the preimage components inside the critical sector of the complementary sector (which are shown here in pale green). 

By analogy with the literature on classical complex dynamics (e.g. \cite{LMM, Sch}), we shall call these sets of pairings the \emph{critical lamination sequence} for the portrait $\mathcal P$ which is realized by an orbit for the connected polynomial sequence $\Pm$. The study of these laminations appears promising as part of the overall program to use symbolic dynamics to describe the behaviour of non-autonomous polynomial sequences and we propose to address this further in a subsequent paper.

\begin{figure}[htbp]

\scalebox{0.35}{\includegraphics{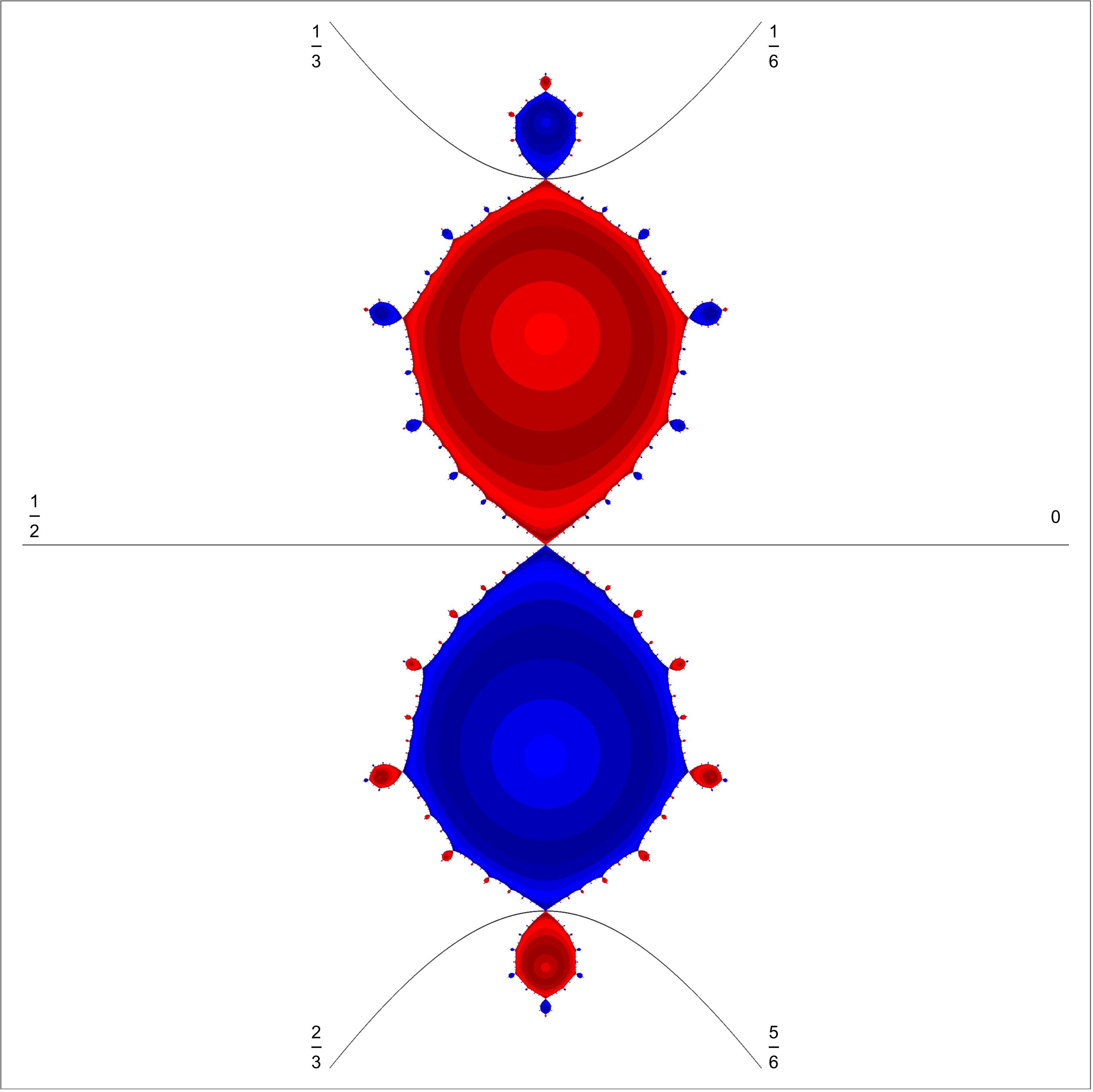}}
\unitlength1cm
\begin{picture}(0.01,0.01)

\end{picture}

\caption{The Juila set for $P(z) = z^2 + \tfrac{3}{2}z$}
\end{figure}

\vspace{.2cm}
{\bf Proof of Proposition 1.1\hspace{.4cm}}The action of $P(z) = z^3 + \tfrac{3}{2}z$ is symmetric under 
the transformations $z \mapsto \overline z$ and $z \mapsto -\overline z$ which correspond to 
reflection in the real and imaginary axes respectively (see Figure 2 below). Additionally, $P$ fixes 
each of these axes. One can check easily that the two critical points at $\pm \tfrac{\im}{\sqrt 2}$ are also fixed and so this polynomial must be hyperbolic with connected Julia set whence all external rays will land on the Julia set and the basin of infinity will be locally connected using either Lemma 2.1 above or the classical results Lemmas 3.2 and 3.3 in \cite{CJY}. 

Recall that the B\"ottcher map $\phi$ which maps the basin of infinity to the outside of the closed unit disc is a limit of the roots of compositions $\sqrt[3^n]{P^{\circ n}}$ using branches of the root functions which are positive on the positive real axis. Now $P$ and all its compositions are odd, symmetric under reflection in both the real and imaginary axes and preserve the positive and negative real axes and it then follows that $\phi$ will inherit these symmetries. 
Additionally, $0$ is a repelling fixed point for $P$ and if $x \ne 0$, then $P^{\circ n}(x)$ is a monotone sequence which converges to $+ \infty$ if $x >0$ and to $-\infty$ if $x < 0$. 
It then follows that $\phi$ must map $(0,\infty)$ and $(-\infty,0)$ to $(1,\infty)$ and $(-\infty, -1)$ rspectively.

For any angle $\theta$, let us denote the (classical) external ray of angle $\theta$ by $R_\theta$. From the above symmetry of the B\"ottcher map, $R_0$ and $R_{1/2}$ land at $0$. 
Again by symmetry, $R_{1/4}$ must lie on the positive imaginary axis and cannot land at $0$ (e.g. because of the superattracting fixed point at $\tfrac{\im}{\sqrt 2}$). 
Hence, if any other external ray of any other angle $\theta$ landed at $0$, we can use symmetry once more to assume without loss of generality that $0 < \theta <\tfrac{1}{4}$. Since 
$\phi$ has derivative $1$ at infinity and rays are not allowed to cross,
such a ray must lie entirely in the first quadrant and $R_0$, $R_{\theta}$ then divide the complex plane into two simply connected regions. Since the complement of the first (closed) quadrant is connected and avoids $R_0 \cup R_\theta \cup \{0\}$, it follows that one and only one of these domains lies entirely in the first quadrant and we shall call this domain $D$. Since $\partial D = R_0 \cup R_\theta \cup \{0\}$, it follows again from the symmetry of $\phi$ and the fact that rays are not allowed to cross that that the closure $\overline D$ can only meet the imaginary axis in a single point, namely $0$. Lastly, since points on the negative real axis have external argument $\tfrac{1}{2}$ and lie outside $D$, it follows from Lemma 1.2 that all external rays with angles in the range $(0, \theta)$ will lie inside $D$.

The ray $R_{\theta/3}$ of angle $\tfrac{\theta}{3}$ thus lies inside $D$ and must accumulate on a preimage of $0$. However, since the other two preimages are $\pm \sqrt{\tfrac{3}{2}}\im$ which lie on the imaginary axis and thus outside $\overline D$, the only allowable preimage is $0$ itself and so this ray must also land at $0$. Continuing in this way, we see that all rays of angles $\tfrac{\theta}{3^n}$, $n \ge 1$ will also land at $0$. However, this is impossible given that this polynomial is hyperbolic which implies that only finitely many rays can land there in view of Lemma 2.2. $\Box$

\vspace{.5cm}
\begin{figure}[htbp]

\scalebox{0.31}{\includegraphics{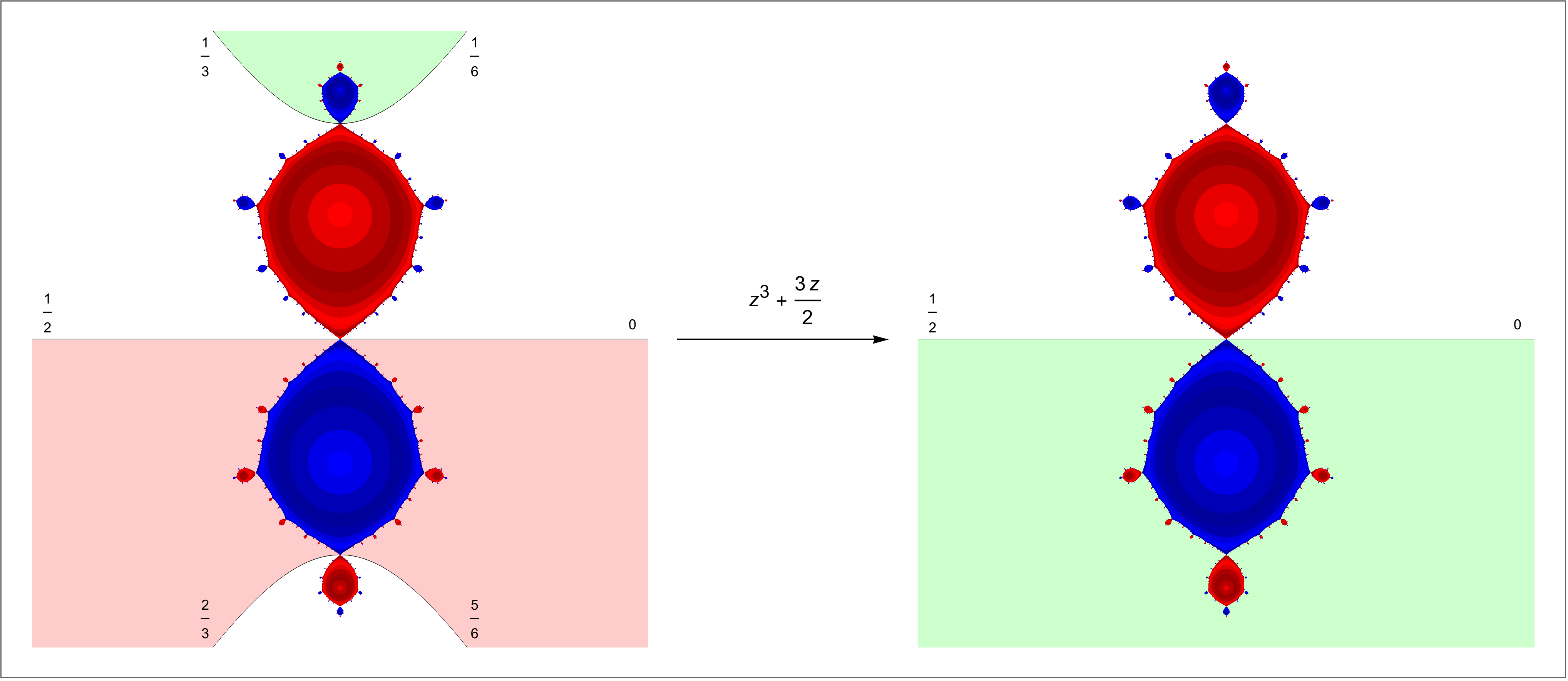}}
\unitlength1cm
\begin{picture}(0.01,0.01)

\end{picture}

\caption{The dynamic portrait at times $0$ and $1$ for the example in Proposition 1.2}
\end{figure}

\vspace{.2cm}
{\bf Proof of Proposition 1.2\hspace{.4cm}}Let $P(z) = z^3 + \tfrac{3}{2}z$ be the same polynomial as in the proof of the last result and let the portrait ${\mathcal P} = \{A_m\}_{m=0}^\infty$ be given by $A_0 = \{\tfrac{1}{6}, \tfrac{1}{3}\}$ and $A_m = \{0, \tfrac{1}{2}\}$ for $m \ge 1$ so that this portrait is invariant under $\theta \mapsto 3\theta \pmod 1$ see Figure 4 above. 

Recall that the external rays $R_0$, $R_{1/2}$ and no others land at $0$.
By the same symmetry of the B\"ottcher map as above, the external rays $R_{1/6}$, $R_{1/3}$ lie entirely in the first and second quadrants respectively and must land on a preimage of $0$. By the same argument as in the proof of Proposition 1.1 above, $R_{1/6}$ cannot land at $0$ and so the only preimage where it can land is $\sqrt{\tfrac{3}{2}}\im$. By symmetry, $R_{1/3}$ must then also land at this point and since only $R_0$, $R_{1/2}$ land at $0$, no rays other than $R_{1/6}$, $R_{1/3}$ can land at $\sqrt{\tfrac{3}{2}}\im$. Similarly the rays $R_{2/3}$, $R_{5/6}$ and no others land at the point $-\sqrt{\tfrac{3}{2}}\im$. 
If we then let $p_0 = \sqrt{\tfrac{3}{2}}\im$ and $p_m = 0$ for $m \ge 1$, we see that this polynomial $P$ with this orbit (hyperbolically) realizes $\mathcal P$. 

By examining their lengths, one sees that the complementary arc $(\tfrac{1}{3}, \tfrac{1}{6})$ for $A_0$ is then critical while the other complementary arc $(\tfrac{1}{6}, \tfrac{1}{3})$ is non-critical. The single critical arc for $A_0$ is then $(\tfrac{1}{3}, \tfrac{1}{6})$ and the image of this arc under $\theta \mapsto 3\theta \pmod 1$ covers the arc $(0, \tfrac{1}{2})$ three times and the other arc $(\tfrac{1}{2},0)$ is covered just twice. Thus the critical value arc for $A_1$ which is associated with $A_0$ must be $(0, \tfrac{1}{2})$. 

By part (1) of Theorem 1.6, since the complementary arc $(\tfrac{1}{6}, \tfrac{1}{3})$ for $A_0$ is non-critical, both critical points for $P$ at $\pm\tfrac{\im}{\sqrt 2}$ must lie in the critical sector associated with the critical arc $(\tfrac{1}{3}, \tfrac{1}{6})$.
Both these critical points are fixed by $P$. Hence one of these corresponding critical values, $-\tfrac{\im}{\sqrt 2}$, is in the lower half plane and thus lies in the sector corresponding to the interval $(\tfrac{1}{2}, 0)$ which is not the one critical value arc for $A_1$ associated with the one critical arc for $A_0$. The result then follows. $\Box$

When dealing with non-autonomous orbit portraits, we also need to grapple with the fact that we no longer have fixed critical arcs and critical value arcs.  The extremely rigid nature of the ray combinatorics for constant sequences enforces a constant angular width on both of these structures at each time.  While we have no \emph{a priori} reason to believe that this is true of non-autonomous portraits, we can establish lower bounds on their sizes, and this is the content of the following two results. The first of these is the main ingredient in allowing us to show that, in the case of constant degree sequences, all non-autonomous orbit portraits which are hyperbolically realizeable are equivalent to periodic portraits.

\vspace{0.2cm}
\begin{lemma}
Let $d, K, M, C, \mu$ be given bounds on degree, coefficients and hyperbolicity. Then there exists $\epsilon>0$ depending only on these constants such that, if 
$\Pm$ is a connected hyperbolic sequence with these bounds and 
$\mathcal{P}$ is any orbit portrait which is associated with an orbit for this sequence then, for each $m \geq 1$, the angular width of any critical value sector at time $m$ for $\mathcal{P}$ is at least $\epsilon$.
 \end{lemma}

\vspace{-.3cm}
\textbf{Proof \hspace{.4cm}}  Suppose for the sake of contradiction that there are constants $d, M,K,C,\mu$ as above for which no such $\epsilon$ exists and that we may find  a sequence sequence of connected hyperbolic sequences $\{ \{P_m^n\}_{m=0}^\infty \}_{n=1}^\infty$ with these constants and a sequence 
of portraits $\{{\mathcal P}^n\}_{n=1}^\infty$ where each ${\mathcal P}^n$ is associated with an orbit for $\Pmn$ and non-negative integers $m_n$ such that the angular width of a critical value sector for ${\mathcal P}^n$ at time $m_n$ converges to $0$ as $n \to \infty$. By truncating the sequences if necessary, we may also assume that each of the times $m_n$ is $1$.

Since all of the polynomial sequences $\Pmn$ above share the same bounds on their degrees, coefficients and hyperbolicity, by Theorem 1.3 of \cite{Com4}, we can find $\delta >0$ such that each of these sequences has postcritical distance $\ge \delta$. 

Let $I^n$ denote these critical value arcs at time $1$ and let $R^{1,n}$ and $R^{2,n}$ denote the rays whose external angles are the endpoints of these critical value arcs. Let $S^n$ denote the corresponding critical value sector enclosed by these rays which contains those points in the basin of infinity at time $1$ with external angles in $I^n$

Let $p^n$ denote the common landing point of $R^{1,n}$, $R^{2,n}$ and for $h >0$ let $x^{1,n}$, $x^{2,n}$ be the points on $R^{1,n}$ and $R^{2,n}$ respectively such that $G_1^n(x_1^n)=G_1^n(x_2^n)=h$ (where of course $G_1^n$ is the corresponding Green's function with pole at $\infty$ for the basin of infinity at time $1$ for the sequence $\{P_m^n\}_{m=0}^\infty$). By Lemma  2.1 we can choose $h$ independently of $n$ such that the Euclidean arclength along $R^{1,n}$ and $R^{2,n}$ from $p_1^n$ to $x^{1,n}$ and $x^{2,n}$ respectively are each less than $\delta/4$. 

Since the iterated Julia sets for all our sequences are connected, the corresponding iterated basins of infinity are simply connected and the level curves on which $G_1^n(z) = h$ are simple closed curves. Now the hyperbolic lengths of the shorter segments of these curves corresponding to external angles in $I^n$ which connect the points $x^{1,n}$ and $x^{2,n}$ must tend to zero as $n$ tends to infinity (we can see this easily using B\"{o}ttcher coordinates and recall that Theorem 3.1 of \cite{CW} gives details concerning the existence and construction of B\"{o}ttcher maps in the general non-autonomous case with connected iterated Julia sets while Theorem 4.7 in \cite{Br} describes the construction in the case of sequences of monic centered quadratic polynomials). By the estimate on the hyperbolic metric at the start of Section 2, we must then have that the Euclidean arc length of this shorter arc of these curves between $x^{1,n}$ and $x^{2,n}$ also tends to $0$ and so we may find $N>0$ such that, for any $n \geq N$,  the (Euclidean) arc length of this curve is less than $\delta/4.$

Since it is a critical value sector, a critical value $c_1^n$ of $P_1^n$ must lie in $S^n$, and the connectedness of the iterated Julia sets implies that it must belong to the corresponding filled Julia set $\mathcal{K}_1^n$ and also to the Fatou set.  The existence of the uniform postcritical distance $\delta$ for the sequence $\{P_m^n\}$ allows us to conclude that the disc ${\mathrm D}(c_1^n,\delta)$ lies in a bounded Fatou component also.

Now any bounded Fatou component must lie in the bounded component of the complement of the Jordan curve corresponding to $G_1^n(z) = h$ which for convenience we shall call $\gamma^n$. The two rays $R^{1,n}$, $R^{2,n}$ and their common landing point $p^n$ give us a crosscut of this Jordan region, which divides it into two smaller Jordan regions. One of these is the Jordan region bounded by the segments of $R^{1,n}$ and $R^{2,n}$ from $p^n$ to $x^{1,n}$ and $x^{2,n}$ and the shorter part of the level curve $\gamma^n$ connecting $x^{1,n}$ and $x^{2,n}$. Let us denote this region by $U^n$ and our claim is that $c_1^n \in U^n$. 

Any curve $\gamma$ from $c_1^n$ to $\infty$ must pass through the level line $\gamma^n$ and it either does so at an external angle belonging to the critical value arc $I^n$ or outside it. Now any point in the critical value sector $S^n$ which is in the basin of infinity can be joined to $c_1^n$ with a curve $\eta$ which avoids $R^{1,n}$ and $R^{2,n}$. We claim that if $\gamma$ first meets $\gamma^n$ outside $S^n$, then $\gamma$ must cross $R^{1,n}$ or $R^{2,n}$ (and at a point where the value of $G^1_n$ is less than or equal to $h$). Otherwise we could combine $\gamma$ and $\eta$ to find a curve which avoided $R_1^n$ and $R^2_n$ but joined a point in the basin of infinity ${\mathcal A}_{\infty,1}^n$ whose external angle lay inside $I^n$ with another in ${\mathcal A}_{\infty,1}^n$ whose external angle lay outside $I^n$. However, this is impossible in view of Lemma 1.2 whence our claim follows.
In either case, whether $\gamma$ crosses $\gamma^n$ inside or outside $S^n$, 
the curve must cross this shorter segment of the Green's line or
these ray segments (or both), whence $c_1^n$ must lie inside $U^n$ in view of the Jordan curve theorem.

By our choice of $h$, $U_n$ has diameter $< \tfrac{\delta}{2}$ and, with the exception of the common landing point of $R^{1,n}$ and $R^{2,n}$, every point on the boundary of $U_n$ is in the corresponding iterated basin of infinity${\mathcal A}_{\infty,1}^n$. On the other hand, from above $c_1^n \in U^n$ and the whole disc ${\mathrm D}(c_1^n, \delta)$ lies in a bounded Fatou component. The result then follows from this contradiction. $\Box$ 

Recall that a critical value arc at time $m+1$ is covered by the image of a critical arc at time $m$ under the $d_{m+1}$-tupling map modulo 1. The following is then an immediate consequence of Part (4) of Theorem 1.6, Lemma 1.1 and the preceding result. 

\vspace{.1cm}
\begin{corollary}
Let $d, K, M,C,\mu$ be given as above. Then if $\mathcal{P}$ is any orbit portrait which can be realized by a connected hyperbolic sequence with these bounds, for each $m \geq 1$, the angular length of each critical value arc at time $m$ for $\mathcal{P}$ is at least $\epsilon$ where $\epsilon >0$ is the number in the statement of the previous result and in particular depends only on the constants $d, K, M,C,\mu$. 

If in addition 
$\mathcal{P}$ is a unicritical orbit portrait, then for each 
$m \ge 0$, we must have that the critical arc of ${\mathcal P}$ at time $m$ has length at least $1 - \tfrac{1 - \epsilon}{d_{m+1}}$. 
\end{corollary}

\vspace{-.1cm}
We are now in a position to establish the equivalence of realizeable non-autonomous portraits with periodic portraits.

\vspace{.2cm}
{\bf Proof of Theorem 1.8 \hspace{.4cm}} Let ${\mathcal P} = \{A_m\}_{m=0}^\infty$ be such an orbit portrait of constant degree which is realized by a connected hyperbolic sequence $\Pm$.  Corollary 2.2 tells us that the critical value arcs must have angular length bounded below by $\epsilon$ where $\epsilon > 0$ is the number in the statement of this corollary. If $d \ge 2$ is the common degree of all the mappings of ${\mathcal P}$, then the iterates of any irrational angle will be dense in the unit circle for the map $\theta \mapsto d\theta \pmod 1$ (note that this is point in the proof where the assumption of constant degree is essential). 

Suppose now we had a critical value arc $I_{m_0}$ of irrational angular length which occurred at some time $m_0 \ge 1$. From above, we can find some time $m_1 > m_0$ such that, at this time, there are two angles in $A_{m_1}$ whose difference is in  
one of the intervals $(\tfrac{k}{d}, \tfrac{k + \epsilon}{d})$ for some $1 \le k \le d-1$ (note that these are intervals on the real line, not the circle). If ${\mathcal P}$ has valence at least $3$, then it follows from the need to preserve cyclic ordering that these two angles will be the endpoints of a single complementary arc. On the other hand, if the valence of ${\mathcal P}$ is two, then again these angles will be the endpoints of one complementary arc as there are simply no other angles in $A_{m_1}$. 

In either case, we have a complementary arc for $A_{m_1}$ whose length is in $(\tfrac{k}{d}, \tfrac{k + \epsilon}{d})$ which (due to the fact that it has length greater than $\tfrac{1}{d}$) must then be a critical arc. The corresponding critical value arc at time $m_1 + 1$ would then have length strictly less than $\epsilon$ which is impossible and this contradiction completes the proof. $\Box$

\vspace{0.2cm}
{\bf Proof of Theorem 1.9 \hspace{.4cm}}Consider ${\mathcal P} = \{A_m\}_{m=0}^\infty$ where $A_m = \{\tfrac{10}{21}, \tfrac{13}{21}, \tfrac{19}{21}\}$ for $m$ even and $A_m = \{\tfrac{17}{21}, \tfrac{20}{21}, \tfrac{5}{21}\}$ for $m$ odd. Note that this is simply the standard portrait $\{\tfrac{1}{7}, \tfrac{2}{7}, \tfrac{4}{7}\}$ shifted by $\tfrac{1}{3}$ at even times and $\tfrac{2}{3}$ at odd times. From this it follows easily that, for each $m \ge 0$,  $A_m$ is mapped bijectively onto $A_{m+1}$ by the doubling map while preserving cyclic ordering and so we do indeed have a formal orbit portrait (of valence $3$). 

There are three possibilities for a smallest possible interval in $\R/\Z$ which contains $\{\tfrac{10}{21}, \tfrac{13}{21}, \tfrac{19}{21}\}$, namely $[\tfrac{10}{21}, \tfrac{19}{21}]$, 
$[\tfrac{13}{21}, \tfrac{10}{21}]$ and $[\tfrac{19}{21}, \tfrac{13}{21}]$. However, each of these contains at least one of $\{\tfrac{17}{21}, \tfrac{20}{21}, \tfrac{5}{21}\}$ so that this cannot give us an unlinked portrait associated with a classical periodic orbit of period $2$ (or even part of a classical portrait of an orbit of higher period). 

To see that this formal portrait can be realized by a hyperbolic non-autonomous sequence of polynomials, let $P_c(z) = z^2 + c$ where $c \approx -.122561 + .744862\im$ is the parameter associated with Douady's rabbit and let $\omega = e^{2\pi \im/3}$. Now let $\phi^0(z)$ and $\phi^1(z)$ be the transformations $z \mapsto \omega z$ and $z \mapsto \omega^2 z$ respectively. Next define a sequence $\Pm$, where $P_m(z) = \phi^1 \circ P_c \circ {(\phi^0)}^{\circ -1}(z) = z^2 + \omega^2c$ for $m$ odd and $P_m(z) = \phi^0 \circ P_c \circ {(\phi^1)}^{\circ -1}(z) = z^2 + \omega c$ for $m$ even. $\Pm$ is then a monic hyperbolic sequence which is conjugate in the non-autonomous sense (e.g. \cite{Com3} Proposition 2.1) to the constant sequence given by $P_c$ using the non-autonomous conjugacy $\phim$ where $\phi_m = \phi^0$ for $m$ even and $\phi_m = \phi^1$ for $m$ odd. 

Let $p$ be the $\beta$-fixed point for $P_c$ and thus the landing point of the rays with angles $\tfrac{1}{7}, \tfrac{2}{7}, \tfrac{4}{7}$. If for each $m \ge 0$ we then let $p_m$ be the image of $p$ under the conjugacy $\phi_m$, this gives an orbit $\{p_m \}_{m=0}^\infty$ for $\Pm$. On the other hand, (e.g. using the uniqueness of the B\"ottcher maps for $\Pm$ as Riemann mappings from the corresponding iterated basins of infinity to $\cbar \setminus \overline \D$ which are tangent to the identity at infinity), we see  that the rays with angles $\tfrac{10}{21}, \tfrac{13}{21}, \tfrac{19}{21}$ (and no others) land at $p_m$ for $m$ even while the rays with angles $\tfrac{17}{21}, \tfrac{20}{21}, \tfrac{5}{21}$ (and no others) land at $p_m$ for $m$ odd. Thus $\Pm$ is a hyperbolic sequence of quadratic polynomials which realizes $\mathcal P$ as desired. $\Box$

\vspace{0.2cm}
{\bf Proof of Theorem 1.10 \hspace{.4cm}} We will construct our sequence from two polynomials $P^0(z)=z^2-1$ and $P^1(z)= z^3$.  Let $\{a_k\}_{k=1}^\infty$ be a binary sequence with entries in $\{0,1\}$ and for $m \ge 1$ define 

\vspace{-.3cm}
\[
P_m(z):=
\begin{cases}
P^0(z) & m \text{ odd}\\
P^{a_k}(z) & m = 2k \:\: \mbox{for some}\:\: k \ge 1.
\end{cases}
\]

The polynomials $P_m$ are all unicritical, and the orbits of the critical points are all bounded, and in fact contained in the set $\{0,-1\}$, so the iterated Julia sets $\mathcal{J}_m$ are all connected.  Moreover, since both these points are either the common critical point $0$ of $P^0$ and $P^1$ or are mapped to $0$ within at most 2 iterations, we may find $\delta>0$ such that, for any $m \ge 0$ and $n\ge 3$, 

\vspace{-.2cm}
\[
Q_{m,m+n}\left(\mathrm{D}(0,\delta)\cup\mathrm{D}(-1,\delta)\right)\subset \mathrm{D}(0,\delta)\cup\mathrm{D}(-1,\delta).
\]

The set $\mathrm{D}(0,\delta)\cup\mathrm{D}(-1,\delta)$ must therefore be contained in some union of (at most two) bounded Fatou components and, since it contains the entire postcritical set $\{0,-1\}$, $\delta$ is a uniform postcritical distance for the polynomial sequence which is independent of the binary sequence $\{a_k\}$ that generated it.  By Theorem 1.3 in \cite{Com4}, the sequence $\{P_m\}$ is hyperbolic with hyperbolicity constants $C$ and $\mu$ which are also independent of the choice of sequence $\{a_k\}$.

By Theorems 1.3 and 3.3 of \cite{Com4}, limit functions on Fatou components must be constant for hyperbolic sequences. On the other hand, $P^0$, $P^1$ both preserve the distance between $-1$, $0$ and so these two points must be in distinct bounded Fatou components for $\Pmn$. Hence for each $m \ge 0$, the real interval $(-1,0)$ must meet the iterated Julia set $\Jm$ and we can find a point $p_m$ 
in $\Jm \cap(-1,0)$. Note also that $P^0$, $P^1$ preserve the interval $(-1, 0)$ and so we can in fact assume that the sequence $\{p_m\}_{m=0}^\infty$ is an orbit for the polynomial sequence $\Pm$.

All of the sequences in question are connected and hyperbolic, so all of the Julia sets are locally connected by \cite{CW} Corollary 3.2. Since the coefficients for every possible sequence are real, it follows by appealing to the resulting symmetry of the B\"ottcher maps in a similar way to in the proof of Proposition 1.1 that, 
for a given angle $\theta \in [0,1)$, the rays with angles $\theta$ and $1-\theta$ will be reflections of each other in the real axis. In particular, the rays for angles $0$, $\tfrac{1}{2}$ must lie on the positive and negative real axes respectively. Also, since $P^0$ and $P^1$ both have real coefficients, they leave the real axis invariant and both map the interval $[-1,1]$ into itself whence this interval belongs to all the iterated filled Julia sets. On the other hand, $P^0$ and $P^1$ are both monotone on each of the intervals $(-\infty, -1)$, $(1, \infty)$ of the complement. It follows from this that all points on the real axis which lie between the landing points of $R_{0,m}$, $R_{1/2,m}$ (which as those points of $\Jm \cap \R$ which are furthest from $0$ must thus have absolute value at least $1$) must have bounded orbits and so the real interval between these landing points is contained in the filled Julia set.

From the above, all rays for angles other than $0$ or $\tfrac{1}{2}$ must then lie entirely in the upper or lower half plane and, given that the B\"ottcher maps are tangent to the identity at $\infty$, it is not hard to see that all rays with angles in $(0,\tfrac{1}{2})$ must lie in the upper half plane while all rays with angles in $(\tfrac{1}{2}, 1)$ must lie in the lower half plane. Thus, if the ray for a given angle $\theta$ lies in one half plane, then the ray for $1-\theta$ will lie in the other half plane. 

By Corollary 2.1 at least one external ray then lands at $p_m \in \Jm$.  By appealing once more to the symmetry of the B\"ottcher maps, we can conclude that for each $m$ we can find at least two rays with angles $\theta_{m}$ and $1-\theta_{m}$ which lie in the upper and lower half  planes respectively and meet at the point $p_m$. Obviously, we can also assume that $\theta_m \in (0,\tfrac{1}{2})$ and the corresponding ray lies in the upper half plane while $1- \theta_m \in (\tfrac{1}{2}, 1)$ and the corresponding ray lies in the lower half plane.  As the points $p_m$ form an orbit (and there are no critical points on the iterated Julia sets), we can thus deduce that these ray pairs belong to an orbit portrait ${\mathcal P}$ for $\Pm$ (note that we do not exclude the possibility of other rays landing on these points).

Now suppose that $\Pm$ and $\Pmt$ are two distinct sequences constructed in this way, and let $m_0$ be the least integer such that $P_{m_0} \neq \tilde{P}_{m_0}$. Note that the first members of both sequences must be $P^0$ so that $m_0 \ge 2$, although we will not make full use of this fact as all we need is that $m_0 \ge 1$ (so that $m_0 -1 \ge 0$). Let $\{p_m\}_{m=0}^\infty$ and $\{\tilde p_m\}_{m=0}^\infty$ respectively be the associated orbits for each of these sequences as above and let $\mathcal P$ and $\tilde {\mathcal P}$ be the associated (formal) portraits for these orbits. 

Now, since $P_{m_0}\neq \tilde{P}_{m_0}$, one of these polynomials must be $P^0$ and one must be $P^1$.  To be definite, let us say that $P_{m_0}=P^0$ and $\tilde{P}_{m_0}=P^1$, and that the rays as above landing at $p_{m_0}$ and $\tilde{p}_{m_0}$ with angles $\theta_{m_0}$ and $\tilde{\theta}_{m_0}$ respectively are in $(0,\tfrac{1}{2})$ and thus lie entirely in the upper half plane.  Now $-1 < p_{m_0}$ and so the critical value $-1$ for $P_{m_0}$ can be joined to the external ray $R_{1/2; m_0}$ by a line segment which avoids $R_{\theta_{m_0};m_0}$ and $R_{1-\theta_{m_0};m_0}$. Thus, in view of Lemma 1.2, we see that the corresponding sector containing those rays with external angles in the interval $(\theta_{m_0},1-\theta_{m_0})$ contains the critical value $-1$ for $P_{m_0}$. Likewise, the sector containing those rays with external angles in the interval $(1-\tilde{\theta}_{m_0},\tilde{\theta}_{m_0})$ contains the critical value $0$ for $\tilde{P}_{m_0}$.  

In a similar way, the sectors associated with their respective polynomial sequences at time $m_0$ with angles in the ranges 
$(1-\theta_{m_0},\theta_{m_0})$ and $(1-\tilde{\theta}_{m_0},\tilde{\theta}_{m_0})$ both contain the critical point $0$ for $P^0 = P_{m_0+1} = \tilde P_{m_0+1}$.
Using Lemma 2.2, we can then choose from among the finitely many rays landing at the point $p_{m_0}$ to ensure that $\theta_{m_0}$ is as large as possible so that $(\theta_{m_0}, 1 - \theta_{m_0})$ is as small as possible. Since $P_{m_0}$ has only one critical point and one critical value, by part (4) of Theorem 1.6, this arc must be the critical value arc for $\mathcal P$ at time $m_0$. However, by part (1) of the same result, it is not the critical arc. In a similar way, we may ensure that $\tilde \theta_0$ is as small as possible in order to ensure that  $(1-\tilde{\theta}_{m_0},\tilde{\theta}_{m_0})$ is the critical value arc and also the critical arc for $\tilde {\mathcal P}$ at time $m_0$.

The arc $(1-\tilde{\theta}_{m_0},\tilde{\theta}_{m_0})$ is the critical value arc for the cubic mapping $\tilde{P}_{m_0}=P^1$ and its preimage under the trebling map will consist of three disjoint preimage arcs of one-third the length each of which covers $(1-\tilde{\theta}_{m_0},\tilde{\theta}_{m_0})$ once under $\theta \mapsto 3\theta \pmod 1$.  The critical arc for $\tilde {\mathcal P}$ at time $m_0 - 1$ must be connected, and so these three preimage arcs will also be connected by two of the preimages of the complementary arc $(\tilde{\theta}_{m_0},1-\tilde{\theta}_{m_0})$.

\begin{figure}[htbp]
\scalebox{0.32}{\includegraphics{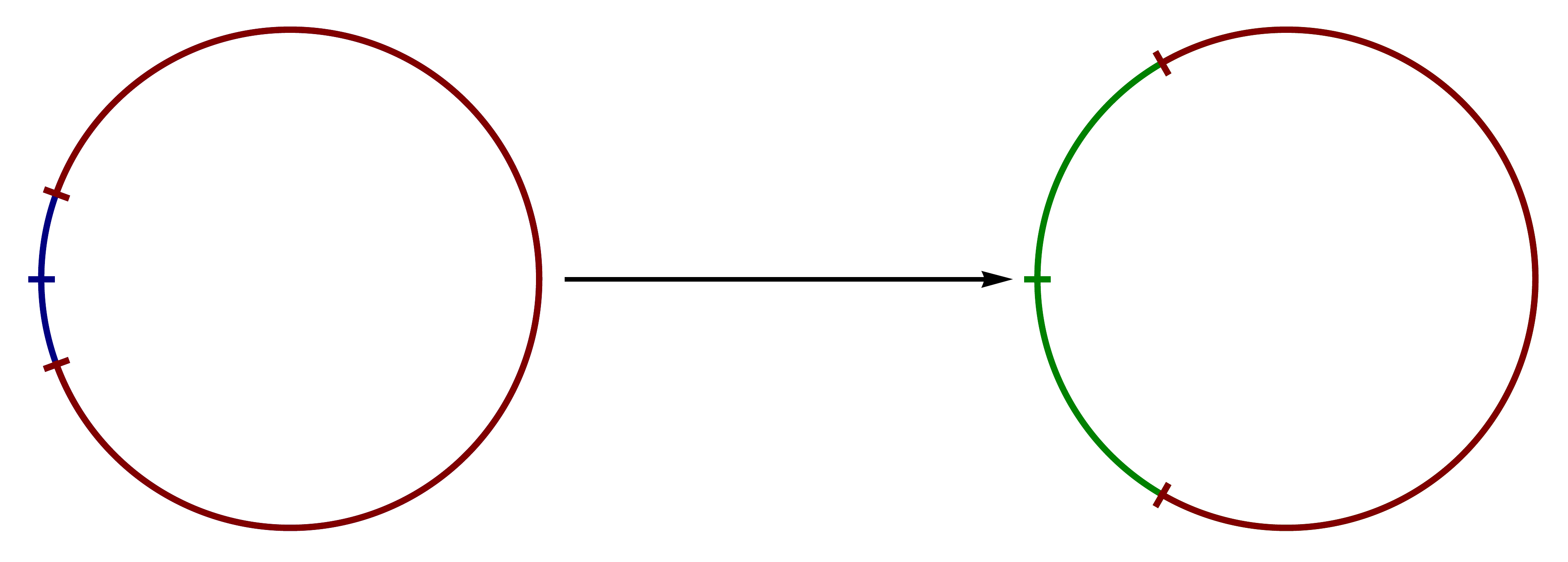}}
\unitlength1cm
\begin{picture}(0.01,0.01)
\put(-7.5,2.45){$\theta \longmapsto 3\theta \pmod 1$}
\put(-10.5,-.35){$\tilde A_{m_0 -1}$}
\put(-2.55,-.35){ $\tilde A_{m_0}$}
\end{picture}

\vspace{.4cm}
\caption{The critical arc at time $m_0 - 1$ in the cubic case.}
\end{figure}
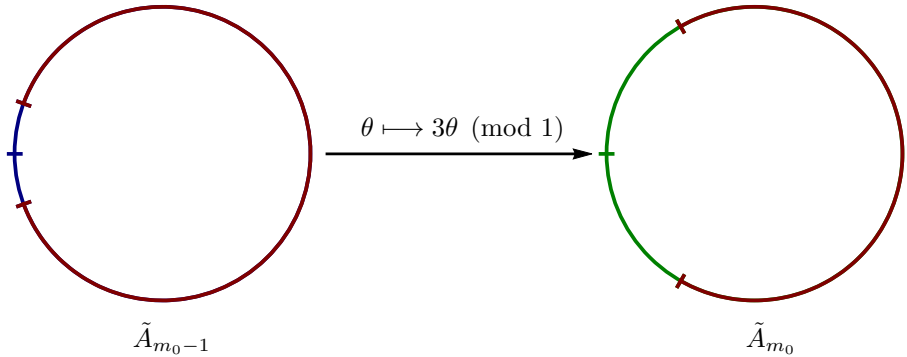

\vspace{.2cm}

The fact that the individual preimages of the critical value arc are arranged symmetrically about the circle informs us that the angular length of the critical arc at time $m_0 - 1$ for $\tilde {\mathcal P}$ will be $\tfrac{2}{3}$ plus two halves of the length of an individual preimage arc.  As it is a critical arc and the next polynomial $P_{m_0 +1}$ in the sequence is $P^0$ which is degree $2$, the arc $(1-\tilde{\theta}_{m_0},\tilde{\theta}_{m_0})$ has length at least $\tfrac{1}{2}$, so its individual preimages under the trebling map will be of length at least $\tfrac{1}{6}$, and this gives us a lower bound of $\tfrac{5}{6}$ on the length of the critical arc at time $m_0 - 1$ for $\tilde {\mathcal P}$ (see Figure 5 above).

The arc $(\theta_{m_0}, 1- \theta_{m_0})$ is a critical value arc but not a critical arc for its portrait $\mathcal P$ at time $m_0$ and so has length strictly less than $\tfrac{1}{2}$.  A similar computation for the preimages of this arc under the doubling map then shows that the critical arc of $\mathcal P$ at time $m_0 -1$ has length less than $\tfrac{3}{4}$ (see Figure 6 below).

\begin{figure}[htbp]
\scalebox{0.32}{\includegraphics{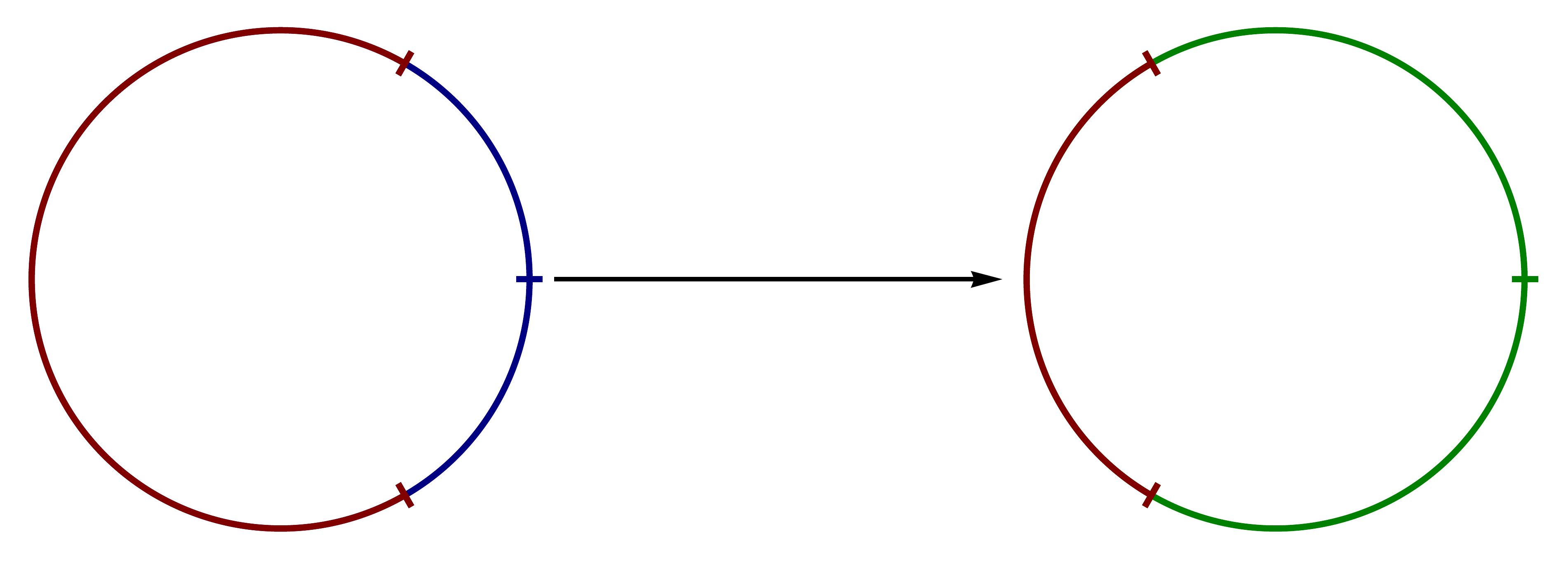}}
\unitlength1cm
\begin{picture}(0.01,0.01)
\put(-7.7,2.45){$\theta \longmapsto 2\theta \pmod 1$}
\put(-10.7,-.35){$A_{m_0 -1}$}
\put(-2.6,-.35){ $A_{m_0}$}
\end{picture}

\vspace{.4cm}
\caption{The critical arc at time $m_0 - 1$ in the quadratic case.}
\end{figure}
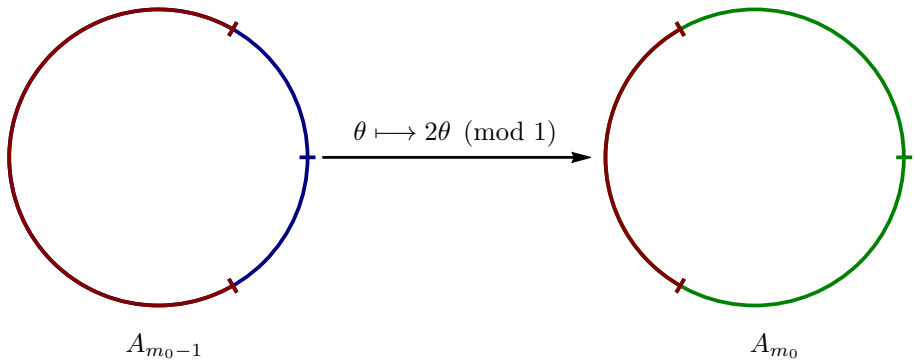

\vspace{.2cm}

Thus the critical arcs at time $m_0 -1$ for ${\mathcal P}$ and $\tilde{\mathcal P}$ are of different length and since these are the only intervals in their respective portraits at this time of length greater than $\tfrac{1}{2}$, each portrait contains a complementary arc at this time of a length which is not present in the other. Since the first $m_0 - 1$ polynomials for both $\Pm$ and $\Pmt$ are the same and in particular have the same degree, it follows that, for any time $0 \le m \le m_{0} -1$, each portrait must contain 
a complementary arc of a length which is not present in the other (since otherwise it would follow that the lengths of all complementary arcs at time $m_0 - 1$ for one of the portraits would occur in the other which we have just shown cannot happen).

Now the set of polynomial sequences $\Pm$ constructed as above is in one-to-one correspondence with the set of binary sequences in $\{0,1\}^{\N}$ which is uncountable. Also, we have just shown that to each distinct such polynomial sequence we can associate a complementary arc at time $0$ whose length does not appear in the complementary arcs at time $0$ for any of the other sequences at this time. On the other hand, by Lemma 2.2, there exists $N \in \N$ which is independent of our choice of sequence such that at most $N$ rays can meet at any point in the Julia set at time $0$ as all these sequences are connected and have uniform degree, coefficient and hyperbolicity bounds. If we could only have rational lengths for the complementary arcs in the portrait at time $0$, this would lead to a contradiction as this would only give us countably many possibilities for the portrait at this time. Hence there must be sequences for which there are complementary arcs of irrational length in the portrait at time $0$ and hence also at all subsequent times. $\Box$

As a final remark, we note that the above example is close to the area of finitely generated polynomial semigroups as studied by Sumi and others. In particular, the existence and properties of the number $\delta$ in the proof shows that the Julia set for the semigroup generated by the two polynomials $P_0 \circ P_0$ and $P_1 \circ P_0$ is hyperbolic in the sense of Definition 1.1 in \cite{Sum1}, which is a stronger statement than merely saying that all the polynomial sequences involved have postcritical distance at least $\delta$. We close with a picture of the Julia set for this semigroup. Note the white areas around $0$ and $-1$ showing how this semigroup is hyperbolic. Since this picture contains the Julia sets for every possible choice, including that for the constant sequence $\{P_0, P_0, P_0, P_0, \ldots \ldots\}$, the interested reader  may like to examine this picture to find the Julia set for this sequence, the well-known basilica. 

\vspace{0.5cm}
\begin{figure}[htbp]
\scalebox{.30}{\includegraphics{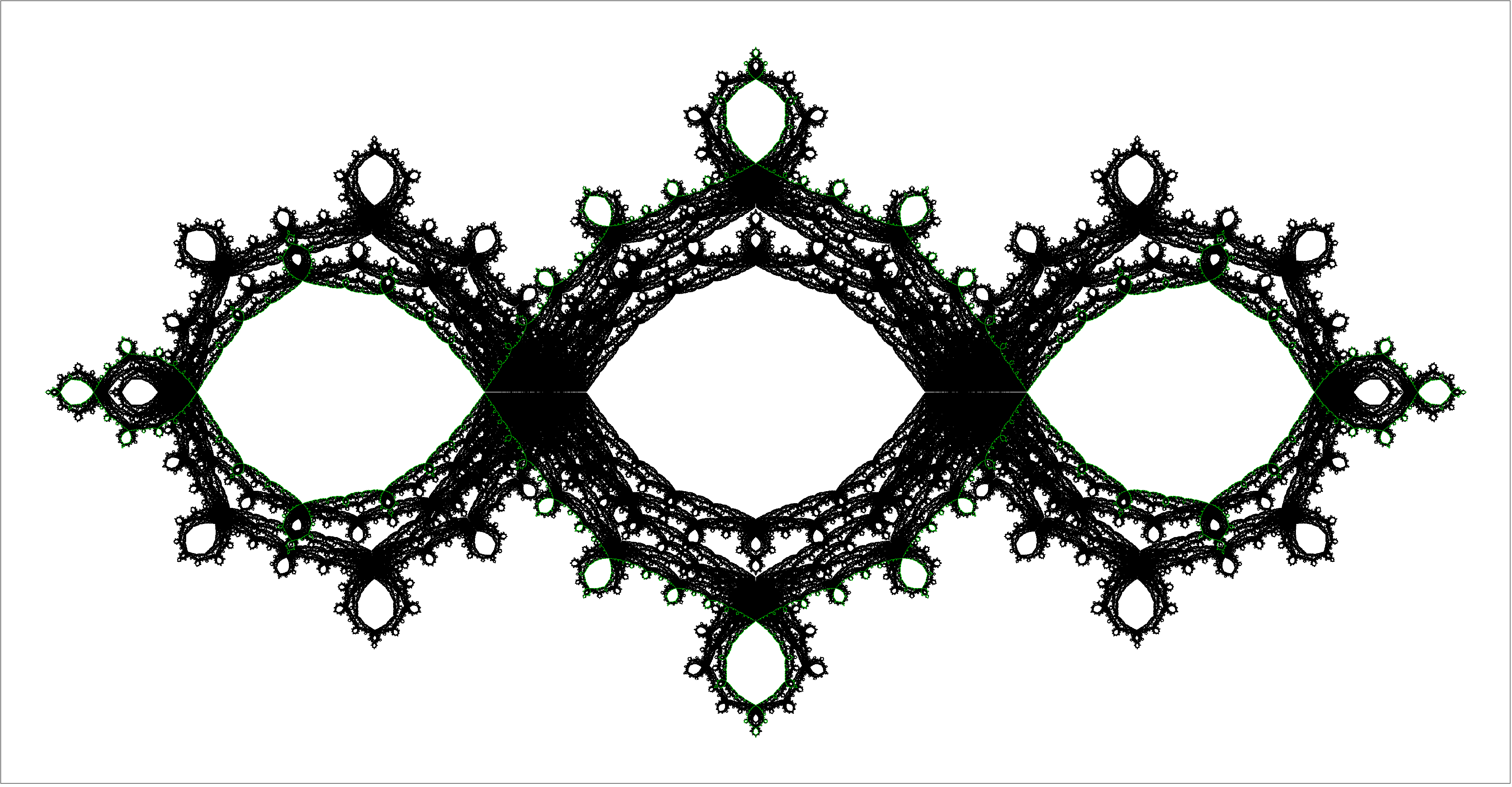}}

\vspace{.4cm}
\caption{A Hyperbolic Semigroup}
\end{figure}

\vspace{.2cm}

\end{document}